\documentclass[12pt]{article}

\usepackage{amsfonts,amssymb,amsmath,amsthm,epsfig,euscript}

\newtheorem{theorem}{Theorem}
\newtheorem{remark}{Remark}
\newtheorem{lemma}[theorem]{Lemma}
\newtheorem{corollary}{Corollary}

\long\def\symbolfootnote[#1]#2{\begingroup
\def\thefootnote{\fnsymbol{footnote}}\footnote[#1]{#2}\endgroup}

\newcommand{\sg}{\sigma}

\def\S{\mathcal{S}}

\title{Classifying Descents According to equivalence mod k}

\author{
Sergey Kitaev \\
\small Institute of Mathematics \\[-0.8ex]
\small Reykjav\'{i}k University \\[-0.8ex]
\small IS-103 Reykjav\'{i}k, Iceland\\[-0.8ex]
\small \texttt{sergey@ru.is}
\and
Jeffrey Remmel \\
\small Department of Mathematics\\[-0.8ex]
\small University of California, San Diego\\[-0.8ex]
\small La Jolla, CA 92093-0112. USA\\[-0.8ex]
\small \texttt{remmel@math.ucsd.edu}
}

\date{\small Submitted: Date 1;  Accepted: Date 2;
 Published: Date 3.\\
\small MR Subject Classifications: 05A15, 05E05}

\begin{document}
\maketitle

\begin{abstract}
In~\cite{kitrem} the authors refine the well-known permutation
statistic ``descent" by fixing parity of (exactly) one of the
descent's numbers. In this paper, we generalize the results
of~\cite{kitrem} by studying descents according to whether the first
or the second element in a descent pair is equivalent to $k$ mod
$k\geq 2$. We provide either an explicit or an inclusion-exclusion type
formula for the distribution of the new statistics. Based on our
results we obtain combinatorial proofs of a number of remarkable identities.
We also provide bijective proofs of some of our
results and state a number of open problems. \\

\end{abstract}

\section{Introduction}
 The {\em descent set}, $Des(\pi)$, of a permutation $\pi=\pi_1\pi_2\cdots\pi_n$ is the set of indices $i$ for which $\pi_i>\pi_{i+1}$. The number of {\em descents} in a permutation $\pi$, denoted by
$des(\pi)$, is a classical permutation statistic. This statistic was first
studied by MacMahon~\cite{macmah} almost a hundred years ago, and
it still plays an important role in the study of permutation statistics.

The {\em Eulerian numbers} $A(n,k)$ count the number of permutations in the
symmetric group $\S_n$ with $k$ descents and they are the
coefficients of the {\em Eulerian polynomials} $A_n(t)$ defined by
$A_n(t)=\sum_{\pi\in\S_n}t^{1+des(\pi)}$. The Eulerian polynomials
satisfy the identity
$$\sum_{k\geq 0}k^nt^k=\frac{A_n(t)}{(1-t)^{n+1}}.$$
For more properties of the Eulerian polynomials see~\cite{comtet}.

In \cite{kitrem}, the authors considered the problem of counting
descents according to the parity of the first or second element of
the descent pair. That is, let $\S_n$ be the set of permutations of
$\{1, \ldots, n\}$, $N = \{0,1,2, \ldots \}$ be the set of natural
numbers, $E = \{0,2,4, \ldots, \}$ be the set of even numbers, $O=
\{1,3,5, \ldots \}$ be the set of odd numbers, and for any statement
$A$, let $\chi(A) =1$ if $A$ is true and $\chi(A)=0$ if $A$ is
false. Then for any $\sg \in \S_n$, define
\begin{itemize}
\item $\overleftarrow{Des}_E(\sg) = \{i: \sg_i > \sg_{i+1} \ \& \
\sg_i \in E\}$ and $\overleftarrow{des}_E(\sg) =
|\overleftarrow{Des}_E(\sg)|$

\item $\overrightarrow{Des}_E(\sg) = \{i: \sg_i > \sg_{i+1} \ \& \
\sg_{i+1} \in E\}$ and $\overrightarrow{des}_E(\sg) =
|\overrightarrow{Des}_E(\sg)|$

\item $\overleftarrow{Des}_O(\sg) = \{i: \sg_i > \sg_{i+1} \ \& \
\sg_i \in O\}$ and $\overleftarrow{des}_O(\sg) =
|\overleftarrow{Des}_O(\sg)|$

\item $\overrightarrow{Des}_O(\sg) = \{i: \sg_i > \sg_{i+1} \ \& \
\sg_{i+1} \in O\}$ and $\overrightarrow{des}_O(\sg) =
|\overrightarrow{Des}_O(\sg)|$
\end{itemize}

Kitaev and Remmel \cite{kitrem} studied the following polynomials:
\begin{enumerate}
\item $R_n(x) = \sum_{\sg \in S_n} x^{\overleftarrow{des}_E(\sg)} = \sum_{k=0}^{n} R_{k,n} x^k$,

\item $P_n(x,z) = \sum_{\sg \in S_n} x^{\overrightarrow{des}_E(\sg)}
z^{\chi(\sg_1 \in E)} = \sum_{k=0}^{n} \sum_{j=0}^1 P_{j,k,n}
z^jx^k$,

\item $M_n(x) = \sum_{\sg \in S_n} x^{\overleftarrow{des}_O(\sg)}
= \sum_{k=0}^{n} M_{k,n} x^k$, and

\item $Q_n(x,z) = \sum_{\sg \in S_n} x^{\overrightarrow{des}_O(\sg)}
z^{\chi(\sg_1 \in O)}=
\sum_{k=0}^{n} \sum_{j=0}^1 Q_{j,k,n} z^jx^k$.
\end{enumerate}

Kitaev and Remmel \cite{kitrem} showed that there are some surprisingly
simple formulas for the coefficients of these polynomials. For example,
they proved
\begin{theorem}\label{thm:KR1}
\begin{eqnarray}
R_{k,2n} &=& \binom{n}{k}^2 (n!)^2, \nonumber \\
R_{k,2n+1}&=& \frac{1}{k+1}{\binom{n}{k}}^2((n+1)!)^2,\nonumber \\
P_{1,k,2n} &=& \binom{n-1}{k} \binom{n}{k+1} (n!)^2,\nonumber \\
P_{0,k,2n} &=& \binom{n-1}{k} \binom{n}{k} (n!)^2,\nonumber \\
P_{0,k,2n+1} &=& (k+1) \binom{n}{k} \binom{n+1}{k+1}(n!)^2 =
(n+1)\binom{n}{k}^2 (n!)^2, \ \mbox{and} \nonumber \\
P_{0,k,2n+1} &=& (n+1)\binom{n}{k}^2(n!)^2. \nonumber
\end{eqnarray}
\end{theorem}

In this paper, we generalize Kitaev and Remmel's results by studying
the problem of counting descents according to whether the first or
the second element in a descent pair is equivalent to $0 \mod k$ for
$k \geq 2$. For any $k > 0$, let $kN = \{0, k, 2k, 3k, \ldots \}$.
Given set $X \subseteq N = \{0, 1, \ldots \}$ and any
$\sg = \sg_1\cdots \sg_n \in \S_n$, we define the following:

\begin{itemize}
\item $\overleftarrow{Des}_X(\sg) = \{i: \sg_i > \sg_{i+1} \ \& \
\sg_i \in X\}$ and $\overleftarrow{des}_X(\sg) =
|\overleftarrow{Des}_X(\sg)|$;

\item $\overrightarrow{Des}_X(\sg) = \{i: \sg_i > \sg_{i+1} \ \& \
\sg_{i+1} \in X\}$ and $\overrightarrow{des}_X(\sg) =
|\overrightarrow{Des}_X(\sg)|$;

\item $A^{(k)}_n(x) = \sum_{\sg \in \S_n}
x^{\overleftarrow{des}_{kN}(\sg)} = \sum_{j=0}^{\lfloor
\frac{n}{k} \rfloor} A^{(k)}_{j,n} x^j$.

\item $B^{(k)}_n(x) = \sum_{\sg \in \S_n}
x^{\overrightarrow{des}_{kN}(\sg)} = \sum_{j=0}^{\lfloor
\frac{n}{k} \rfloor} B^{(k)}_{j,n} x^j$.

\item $B^{(k)}_n(x,z) = \sum_{\sg \in \S_n}
x^{\overrightarrow{des}_{kN}(\sg)}z^{\chi(\sg_1 \in kN)} =
\sum_{j=0}^{\lfloor \frac{n}{k} \rfloor} \sum_{i=0}^1
B^{(k)}_{i,j,n} z^ix^j$.

\end{itemize}

\begin{remark} Note that setting $k=1$ gives us {\rm(}usual{\rm)} descents providing
$A^{(1)}_n(x)=B^{(1)}_n(x)=A_n(x)$, whereas setting $k=2$ gives
$\overleftarrow{Des}_E(\sg)$ and $\overrightarrow{Des}_E(\sg)$
studied in~$\cite{kitrem}$.
\end{remark}

The goal of this paper is to derive closed formulas for the coefficients of
these polynomials. When $k > 2$, our  formulas are considerably more complicated than the formulas in the $k=2$ case. In fact, in most cases, we can derive two distinct formulas for the coefficients of these polynomials. We shall see that there are simple
recursions for the coefficients of the
polynomials $A^{(k)}_{kn+j}(x)$, $B^{(k)}_{kn+j}(x)$, and  $B^{(k)}_{kn+j}(x,z)$ for $0 \leq j \leq k-1$. In fact, we can derive
two different formulas for the coefficients of our polynomials
by iterating the recursions
starting with the constant term  and by
iterating the recursions starting with highest coefficient.
For example, we shall prove the following theorem.

\begin{theorem}\label{thm:KR2}
For all $k \geq 2$, $n \geq 0$, and $0 \leq j \leq k-1$,
\begin{eqnarray}
&&A^{(k)}_{s,kn+j} = \nonumber \\
&&((k-1)n+j)!  \sum_{r=0}^s (-1)^{s-r}\binom{(k-1)n+j+r}{r}
\binom{kn+j+1}{s-r} \prod_{i=0}^{n-1} (r+1+j +(k-1)i) = \nonumber \\
&& ((k-1)n+j)!  \sum_{r=0}^{n-s} (-1)^{n-s-r}\binom{(k-1)n+j+r}{r}
\binom{kn+j+1}{n-s-r} \prod_{i=1}^{n} (r+(k-1)i) \nonumber
\end{eqnarray}
\end{theorem}

What is remarkable about these two different formulas for
$A^{(k)}_{kn+j}$ is that they lead to a number of identities that
are interesting in the own right. For example, it follows from
Theorem \ref{thm:KR2} that for all $k \geq 2$, $n \geq 0$, and $0
\leq j \leq k-1$,
\begin{eqnarray}\label{special1}
&&\sum_{r=0}^s (-1)^{s-r}\binom{(k-1)n+j+r}{r}
\binom{kn+j+1}{s-r} \prod_{i=0}^{n-1} (r+1+j +(k-1)i) = \nonumber \\
&&\sum_{r=0}^{n-s} (-1)^{n-s-r}\binom{(k-1)n+j+r}{r}
\binom{kn+j+1}{n-s-r} \prod_{i=1}^{n} (r+(k-1)i) \nonumber
\end{eqnarray}
Even in the case $k=2$, we get some remarkable identities. For
example, it follows from Theorems \ref{thm:KR1} and \ref{thm:KR2}
that  for all $n \geq s$,
\begin{eqnarray}\label{special2}
\binom{n}{s}^2(n!) &=& \sum_{r=0}^s (-1)^{s-r}\binom{n+r}{r}
\binom{2n+1}{s-r} \prod_{i=0}^{n-1} (r+1+i) \nonumber \\
&=& \sum_{r=0}^{n-s} (-1)^{n-s-r}\binom{n+r}{r} \binom{2n+1}{n-s-r}
\prod_{i=1}^{n} (r+i)
\end{eqnarray}
It turns out that both of these identities can be derived by using
certain hypergeometric series identities. For example, we will show
how (\ref{special2}) can be derived from Saalcsh\"{u}tz's identity.
Jim Haglund \cite{haglund} suggested that (\ref{special1}) should
follow from Gasper's transformation \cite{Gasper} of hypergeometric
series of Karlsson-Minton type. This is indeed the case but we will
not include such a derivation in this paper since (\ref{special1})
is a special case of wider class of identities that arise by
studying the problem of enumerating permutations according to the
number of pattern matches where the equivalence classes of the
elements modulo $k$ for $k \geq 2$ are taken into account, see
\cite{liese}. A general derivation of this wider class of identities
from the Gasper's transformation of hypergeometric series of
Karlsson-Minton type will appear in a subsequent paper
\cite{liese-rem}.

Given any permutation $\sg = \sg_1\cdots \sg_n \in \S_n$, we label
the possible positions of where we can insert $n+1$ to get a
permutation in $\S_{n+1}$ from left to right with 0 to $n$, i.e.,
inserting $n+1$ in position 0 means that we insert $n+1$ at the
start of $\sg$ and for $i \geq 1$, inserting $n+1$ in position $i$
means we insert $n+1$ immediately after $\sg_i$. In such a
situation, we let $\sg^{(i)}$ denote the permutation of $\S_{n+1}$
that results by inserting $n+1$ in position $i$.

Let $\sg^c = (n+1 -\sg_1)(n+1 -\sg_2) \cdots (n+1 -\sg_{n})$
denote the {\em complement} of $\sg$. Clearly,
if $n$ is odd, then, for all $i$,  $\sg_i$ and
$n+1 -\sg_i$ have the same parity,
whereas they have opposite parity if $n$ is even. However, if $k \geq 3$,
then complementation does not preserve equivalences classes mod $k$.
The {\em
reverse} of $\sg$ is the permutation
$\sg^r=\sg_{n}\sg_{n-1}\cdots\sg_1$.

The outline of this paper is as follows. In section~\ref{section2},
we shall give explicit formulas for the coefficients
$A^{(k)}_{0,kn+j}$ and $A^{(k)}_{n,kn+j}$ for all $k \geq 2$, $n
\geq 0$, and $j \in \{0, \ldots,k-1\}$. Then we shall develop a set
of recursions for the coefficients $A^{(k)}_{s,kn+j}$ and use these
recursions to derive our two different formulas for the coefficients
of $A^{(k)}_{kn+j}(x)$ for all $n \geq 0$ and $j \in \{0,
\ldots,k-1\}$. In section~\ref{section3}, we shall give explicit
formulas for the coefficients $B^{(k)}_{0,kn+j}$,
$B^{(k)}_{0,0,kn+j}$,$B^{(k)}_{1,0,kn+j}$, $B^{(k)}_{n,kn+j}$,
$B^{(k)}_{0,n,kn+j}$, and $B^{(k)}_{1,n,kn+j}$ for all $k \geq 2$,
$n \geq 0$, and $j \in \{0, \ldots,k-1\}$. Then we shall develop a
set of recursions for the coefficients $B^{(k)}_{i,s,kn+j}$ and use
these recursions to derive inclusion-exclusion type formulas for the
coefficients of $B^{(k)}_{s,kn+j}$, $B^{(k)}_{0,s,kn+j}$, and
$B^{(k)}_{1,s,kn+j}$ for all $n \geq 0$ and $j \in \{0,
\ldots,k-1\}$. Based on such formulas, we shall derive a number of
remarkable identities (See Theorem~\ref{thm:Bid1}). In
section~\ref{section4}, we shall consider some natural bijective
questions that arise from our results. Finally, in
section~\ref{section5}, we shall discuss a number of open questions.

\section{Properties of $A_n^{(k)}(x)$}\label{section2}

In this section, we shall study the properties of the polynomials
$A_n^{(k)}(x)$. For example, here are some examples of the polynomials
$A_n^{(3)}(x)$.

\begin{description}

\item $A^{(3)}_1(x) =1$.

\item $A^{(3)}_2(x) =2$.

\item $A^{(3)}_3(x) =2+4x$.

\item $A^{(3)}_4(x) =12+12x$.

\item $A^{(3)}_5(x) =72+48x$.

\item $A^{(3)}_6(x) =72+456x+192x^2$.

\item $A^{(3)}_7(x) =960+3120x+960x^2$.

\item $A^{(3)}_8(x) =10800+ 23760x+5760x^2$.

\item $A^{(3)}_9(x) =10800 + 133920x + 183600x^2 + 34560x^3$.

\item $A^{(3)}_{10}(x) =241920 + 1572480x + 1572480x^2 + 241920 x^3$.

\item $A^{(3)}_{11}(x) =4233600 + 18869760x + 14878080x^2 +
1935360x^3$.

\item $A^{(3)}_{12}(x) =4233600 + 84309120x + 233331840x^2 + 141644160 x^3
15482880x^4$.

\item $A^{(3)}_{13}(x) = 139345920 + 1478373120x + 2991582720x^2 +
1478373120x^3 + 139345920x^4$.

\item $A^{(3)}_{14}(x) =3429216000 + 25202016000x + 40334112000x^2 +
16819488000x^3 + 1393459200 x^4 $.

\item $A^{(3)}_{15}(x) =3429216000 + 98413056000x + 448628544000x^2 +
551287296000x^3 + 191981664000x^4 + 1393459200x^5$.

\end{description}

By Theorem \ref{thm:KR2}, we know that $A^{(3)}_{s,3n+j}$
is divisible by $(2n+j)!$. However, if we consider
$A^{(3)}_{4,15}/(10!) = 191981664000/(10!) = 52905$, then one can
check that the prime factorization of 52905 is $3\cdot 5 \cdot 3527$.
Thus the prime $3527$ divides $A^{(3)}_{4,15}$ so that we can not expect
that we will get formulas for $A^{(3)}_{s,3n+j}$ that are as
 simple as the formulas that appear in Theorem \ref{thm:KR1} for the
polynomials $A^{(2)}_{s,2n+j}$.

For the rest of this paper, we shall assume that $k \geq 2$.

For $j = 1, \ldots, k-1$, let $\Delta_{kn+j}$ be the operator
which sends $x^s$ to $s x^{s-1} + (kn+j-s)x^s$ and $\Gamma_{kn+k}$
be the operator that sends $x^s$ to $(s+1)x^s +(kn+k-1-s)
x^{s+1}$. Then we have the following.

\begin{theorem}\label{thm:1}
The polynomials $\{A^{(k)}_n(x)\}_{n \geq 1}$ satisfy the following recursions.

\begin{enumerate}
\item $A^{(k)}_1(x) = 1$, \item For $j=1, \ldots, k-1$,
$A^{(k)}_{kn+j}(x) = \Delta_{kn+j}(A^{(k)}_{kn+j-1}(x))$ for $n
\geq 0$, and \item $A^{(k)}_{kn+k}(x) =
\Gamma_{kn+k}(A^{(k)}_{kn+k-1}(x))$ for $n \geq 1$.
\end{enumerate}
\end{theorem}
\begin{proof}
Part (1) is trivial.

For part (2), fix $j$ such that $1 \leq j \leq k-1$. Now suppose
$\sg = \sg_1 \cdots \sg_{kn+j-1} \in \S_{kn+j-1}$ and
$\overleftarrow{des}_{kN}(\sg) =s$. It is then easy to see that if
we insert $kn+j$ in position $i$ where $i \in
\overleftarrow{Des}_{kN}(\sg)$, then
$\overleftarrow{des}_E(\sg^{(i)}) =s-1$. However, if we insert
$kn+j$ in position $i$ where $i \notin
\overleftarrow{Des}_{kN}(\sg)$, then
$\overleftarrow{des}_{kN}(\sg^{(i)}) =s$.  Thus $\{\sg^{(i)}: i =
0, \ldots, kn+j-1\}$ gives a contribution of $sx^{s-1} + (kn+j
-s)x^s$ to $A^{(k)}_{kn+j}(x)$.

For part (3), suppose $\sg = \sg_1 \cdots \sg_{kn+k-1} \in
\S_{kn+k-1}$ and $\overleftarrow{des}_{kN}(\sg) =s$. It is then
easy to see that if we insert $kn+k$ in position $i$ where $i \in
\overleftarrow{Des}_E(\sg)$ or $i=kn+k-1$, then
$\overleftarrow{des}_E(\sg^{(i)}) =s$. Similarly if we insert
$kn+k$ in position $i$ where $i \notin
\overleftarrow{Des}_{kN}(\sg) \cup \{kn+k-1\}$, then
$\overleftarrow{des}_{kn}(\sg^{(i)}) =s+1$. Thus $\{\sg^{(i)}: i =
0, \ldots, kn+k-1\}$ gives a contribution of $(s+1)x^{s} + (kn+k
-(s+1))x^{s+1}$ to $A^{(k)}_{kn+k}(x)$.
\end{proof}

Note that we can rewrite Theorem~\ref{thm:KR1} as saying that
\begin{enumerate}
\item $A^{(k)}_1(x) = 1$, \item For $j=1, \ldots, k-1$,
$A^{(k)}_{kn+j}(x) = (1-x)\frac{d}{dx}(A^{(k)}_{kn+j-1}(x)) +
(kn+j)A^{(k)}_{kn+j-1}(x)$ for $n \geq 0$, and \item
$A^{(k)}_{kn+k}(x) = (x -x^2)\frac{d}{dx}(A^{(k)}_{kn+k-1}(x))+
(1+x(kn+k-1))A^{(k)}_{kn+k-1}(x)$ for $n \geq 1$.
\end{enumerate}

There are simple formula for lowest and the highest coefficients
in the polynomials $A^{(k)}_{kn+j}$ and we can give direct combinatorial
proofs of such formulas. That is, we can prove the following.

\begin{theorem}\label{thm:2}  We have
\begin{itemize}
\item[{\rm(}a{\rm)}] $A^{(k)}_{0,kn+j} = ((k-1)n+j)!
\prod_{i=0}^{n-1}(j+1 + i(k-1))$ for $0\leq j\leq k-1$;
\item[{\rm(}b{\rm)}] $A^{(k)}_{n,kn+j} = (n(k-1)+j)! (k-1)^n n!$
for $0\leq j\leq k-1$.
\end{itemize}
\end{theorem}
\begin{proof}
It is easy to see that both (a) and (b) hold for $n=1$.

To prove (a) fix $j$ such that $0 \leq j \leq k-1$ and suppose
that $\sg = \sg_1 \cdots \sg_{kn+j}$ is such that
$\overleftarrow{des}_{kN}(\sg) =0$. Then we can factor any such
permutation into blocks by reading the permutation from left to
right and cutting after each  number which is not equal to $0 \mod
k$. For example if $j=0$, $k=3$, and $\sg =
11~1~2~4~5~3~6~7~8~9~10~12$, then the blocks of $\sg$ would be
$11$, $1$, $2$,  $4$, $5$, $3~6~7$,  $8$, $9~10$, $12$.

There may be a block of numbers which are equivalent to $0 \mod k$
at the end which are arranged in increasing order. We call this
final block the $\infty$-th block. Every other block must end with
a number $sk+i$ where $0 \leq s \leq n-1$ and $1 \leq i \leq k-1$
which can be preceded by any subset of numbers which are
equivalent to $0 \mod k$ and which are less than $sk+i$ arranged
in increasing order. We call such a block the $(sk+i)$-th block.
It is then easy to see that there are $\prod_{i=0}^{n-1}
(j+1+i(k-1))$ ways to put the numbers $k, 2k, \ldots, nk$ into the
blocks.  That is, $kn$ must go in either the blocks $kn+1, \ldots,
kn+j$ or the $\infty$-th block so that there are $1+j$ choices for
the block in which to place $kn$.  Then $k(n-1)$ can either go in the
blocks $k(n-1)+1, k(n-1)+2, \ldots, k(n-1)+k-1, kn+1, \ldots, kn+j$
or the $\infty$-block  so that there are $j+1 +(k-1)$ choices for the
block that contains $k(n-1) $. More generally, $k(n-i)$ can go in
any blocks $k(n-i)+1, \ldots, k(n-i)+k-1, k(n-i+1)+1, \ldots,
k(n-i+1)+k-1, \ldots , k(n-1)+1, \ldots, k(n-1)+k-1, kn+1, \ldots,
kn+j$ or the $\infty$-block so that are $j+1 +i(k-1)$ choices for
the block that contains $k(n-i)$.  Once we have arranged the
numbers which are equivalent to $0 \mod k$ into blocks, it is easy
to see that we can arrange blocks $sk+i$ where $0 \leq s \leq n-1$
and $1 \leq i \leq k-1$ plus the blocks $kn+1, \ldots, kn+j$ in
any order and still get a permutation $\sg$ with
$\overleftarrow{des}_{kN}(\sg)= 0$. It thus follows that there are
$((k-1)n+j)! \prod_{i=0}^{n-1} (j+1 + i(k-1))$ such permutations.

To prove (b) fix $j$ such that $0 \leq j \leq k-1$ and suppose
that $\sg = \sg_1 \cdots \sg_{kn+j}$ is such that
$\overleftarrow{des}_{kN}(\sg) =n$. Then, as above, we can factor
any such permutation into blocks by reading the permutation from
left to right and cutting after each number which is not equal to
0 mod $k$. One can see that, unlike the case where
$\overleftarrow{des}_{kN}(\sg) =0$, there can be no numbers which
are equivalent to 0 mod $k$ at the end (otherwise not all of such
numbers will be involved in a descent). Thus the $\infty$-th block
must be empty. Similarly, it is easy to see that the blocks $kn+1,
\ldots, kn+j$ must be singletons. Next if $0 \leq s \leq n-1$ and $1
\leq j \leq k-1$ and there are numbers which are equivalent to $0
\mod k$ in the $(sk+j)$-th block, i.e. the block that ends with
$sk+j$, then those numbers must all be greater than $sk+j$ and
they must be arranged in decreasing order.

It is then easy to see that there are $(k-1)^n (n!)$ ways to put
the numbers $k, 2k, \ldots, nk$ into blocks.  That is, $kn$ may go
in any of the blocks $sk+j$ where $0 \leq s \leq n-1$ and $1 \leq
j \leq k-1$  so that there are $(k-1)n$ choices for the block that
contains $kn$. Then $k(n-1)$ can go in any of the blocks $sk+j$
where $0 \leq s \leq n-2$ and $1 \leq j \leq k-1$  so that there
are $(k-1)(n-1)$ choices for the block that contains $k(n-1)$,
etc. After we have partitioned the  numbers which are equivalent
to $0 \mod k$ into their respective blocks, we must arrange the
numbers which are equivalent to $0 \mod k$ in each block in
decreasing order so that there are a total $(k-1)^n n!$ ways to
partition the numbers which are equivalent to $0 \mod k$ into the
blocks. Once we have arranged the numbers which are equivalent to
$0 \mod k$ into blocks, it is easy to see that we can arrange
blocks $sk+j$ where $0 \leq s \leq n-1$ and $1 \leq j \leq k-1$
plus the blocks $kn+1, \ldots, kn+j$ in any order an still get a
permutation $\sg$ with $\overleftarrow{des}_{kN}(\sg)= n$. It thus
follows that there are $(n(k-1)+j)! (k-1)^n n!$ such permutations.
\end{proof}

It is easy to see from Theorem \ref{thm:1} that we have two
following recursions for the coefficients $A^{(k)}_{s,n}$.

For $1 \leq j \leq k-1$,
\begin{equation}\label{kn+j-1:kn+jrec}
A^{(k)}_{s,kn+j} =  (kn+j-s)A^{(k)}_{s,kn+j-1} +(s+1)A^{(k)}_{s+1,kn+j-1}
\end{equation}
and
\begin{equation}\label{kn+k-1:kn+krec}
A^{(k)}_{s,kn+k} =  (1+s)A^{(k)}_{s,kn+k-1} +(kn+k-s)A^{(k)}_{s-1,kn+k-1}.
\end{equation}

The following theorem provides an inclusion-exclusion type formula
for $A^{(k)}_{s,kn+j}$ which can be obtained by interating the
recursions \ref{kn+j-1:kn+jrec} and \ref{kn+k-1:kn+krec} starting
with our formulas for $A^{(k)}_{0,kn+j}$.

\begin{theorem}\label{thm:5}

For all $0 \leq j \leq k-1$ and all $n \geq 0$, we
have
\begin{eqnarray*}
&&A^{(k)}_{s,kn+j} = \\
&&((k-1)n+j)! \left[ \sum_{r=0}^s (-1)^{s-r}\binom{(k-1)n+j+r}{r}
\binom{kn+j+1}{s-r} \prod_{i=0}^{n-1} (r+1+j +(k-1)i)\right]. \nonumber
\end{eqnarray*}
\end{theorem}
\begin{proof}
We shall prove this formula by induction on $s$. Note that
Theorem~\ref{thm:2} shows that our formula for $A^{(k)}_{s,kn+j}$ holds
when $s = 0$ for
all $n \geq 0$ and $0 \leq j \leq k-1$.

Now assume by induction that our formula for $A^{(k)}_{s,kn+j}$ is
true for  all $n \geq 0$ and $0 \leq j \leq k-1$. Then we shall
prove that it holds for $A^{(k)}_{s+1,kn+j}$ for all $n \geq 0$ and
$0 \leq j \leq k-1$. Note that by recursion  (\ref{kn+j-1:kn+jrec}),
we have for $1 \leq j \leq k-1$,
\begin{equation*}
(s+1)A^{(k)}_{s+1,kn+j-1} = A^{(k)}_{s,kn+j} - (kn+j-s)A^{(k)}_{s,kn+j-1}.
\end{equation*}
Thus
\begin{eqnarray}
&&(s+1)A^{(k)}_{s+1,kn+j-1} = \nonumber \\
&&((k-1)n+j)! \left[ \sum_{r=0}^s (-1)^{s-r} \binom{(k-1)n+j+r}{r}
\binom{kn+j+1}{s-r} \prod_{i=0}^{n-1} (r+1+j +(k-1)i)\right] \nonumber \\
&&-(kn+j-s) ((k-1)n+j-1)! \times \nonumber \\
&& \ \ \ \ \  \ \ \left[ \sum_{r=0}^s (-1)^{s-r}
\binom{(k-1)n+j-1+r}{r} \binom{kn+j}{s-r} \prod_{i=0}^{n-1} (r+ j
+(k-1)i)\right].\nonumber
\end{eqnarray}
It follows that
\begin{eqnarray}\label{genform4}
&&\frac{A^{(k)}_{s+1,kn+j-1}}{((k-1)n+j-1)!} = \\
&&\frac{((k-1)n+j)}{s+1} \left[ \sum_{r=0}^s (-1)^{s-r} \binom{(k-1)n+j+r}{r}
\binom{kn+j+1}{s-r} \prod_{i=0}^{n-1} (r+1+j +(k-1)i)\right] \nonumber \\
&&-\frac{(kn+j-s)}{s+1} \left[ \sum_{r=0}^s (-1)^{s-r}
\binom{(k-1)n+j-1+r}{r} \binom{kn+j}{s-r} \prod_{i=0}^{n-1} (r+ j
+(k-1)i)\right].\nonumber
\end{eqnarray}
We can divide the terms on the RHS of (\ref{genform4}) into the three
parts.
The $r=s$ term from the first summand on the RHS of (\ref{genform4}) gives
\begin{eqnarray}\label{genform5}
&&\frac{((k-1)n+j)}{s+1} \binom{(k-1)n+j+s}{s}
\prod_{i=0}^{n-1} (s+1+j +(k-1)i)=\\
&& \binom{(k-1)n+j+s}{s+1}
\prod_{i=0}^{n-1} (s+1+j +(k-1)i) =  \nonumber \\
&&\binom{(k-1)n+(j-1)+s+1}{s+1} \prod_{i=0}^{n-1} ((s+1)+ 1+ (j-1)
+(k-1)i).\nonumber
\end{eqnarray}

The $r=0$ term from the second  summand on the RHS of (\ref{genform4}) gives
\begin{eqnarray}\label{genform6}
&&(-1)^{s+1}\frac{(kn+j-s)}{s+1} \binom{kn+j}{s}
\prod_{i=0}^{n-1} (j +(k-1)i)=\\
&& (-1)^{s+1}\binom{kn+j}{s+1}
\prod_{i=0}^{n-1} (j +(k-1)i) =  \nonumber \\
&&(-1)^{s+1} \binom{kn+(j-1)+1}{s+1} \prod_{i=0}^{n-1} (1+ (j-1)
+(k-1)i).\nonumber
\end{eqnarray}

Finally we can organize the remaining terms of each summand
according to the factor $\prod_{i=0}^{n-1} (r+j+(k-1)i)$ to get

\begin{eqnarray}\label{genform7}
&&\sum_{r=1}^s (-1)^{s+1-r} \prod_{i=0}^{n-1} (r+j+(k-1)i) \times
\left[\frac{kn+j-s} {s+1} \binom{(k-1)n+j-1+r}{r}
\binom{kn+j}{s-r} + \nonumber \right. \\
&& \left. \frac{(k-1)n+j}{s+1} \binom{(k-1)n+j+r-1}{r-1}
\binom{kn+j+1}{s-(r-1)}\right] =
\end{eqnarray}

\begin{eqnarray} && \sum_{r=1}^s (-1)^{s+1-r}
\frac{((k-1)n+j+r-1)\downarrow_r}{(r-1)!}
\frac{(kn+j)\downarrow_{s-r}}{(s-r)!} \frac{1}{s+1}
\prod_{i=0}^{n-1} (r+j+(k-1)i) \times \nonumber \\
&& \left[\frac{kn+j-s}{r} + \frac{kn+j+1}{s+1-r}\right] = \nonumber
\end{eqnarray}

\begin{eqnarray}
&&\sum_{r=1}^s (-1)^{s+1-r}
\frac{((k-1)n+j+r-1)\downarrow_r}{(r-1)!}
\frac{(kn+j)\downarrow_{s-r}}{(s-r)!} \frac{1}{s+1}
\prod_{i=0}^{n-1} (r+j+(k-1)i) \times \nonumber \\
&&\frac{(s+1)(kn+j-s+r)}{r(s+1-r)} = \nonumber
\end{eqnarray}

\begin{eqnarray}
&&\sum_{r=1}^s (-1)^{s+1-r} \binom{(k-1)n+j+r-1}{r}
\binom{kn+j}{s+1-r}
\prod_{i=0}^{n-1} (r+j+(k-1)i) = \nonumber \\
&&\sum_{r=1}^s (-1)^{s+1-r} \binom{(k-1)n+j-1 +r}{r}
\binom{kn+(j-1)+1}{s+1-r} \prod_{i=0}^{n-1} (r+1+(j-1)+(k-1)i),
\nonumber
\end{eqnarray}
where $x\downarrow_{r}=x(x-1)\cdots (x-r+1)$. Thus if we combine
(\ref{genform5}), (\ref{genform6}), and (\ref{genform7}), we get
\begin{eqnarray*}
&&\frac{A^{(k)}_{s+1,kn+j-1}}{((k-1)n+j-1)!} = \\
&&\sum_{r=0}^{s+1} (-1)^{s+1-r} \binom{(k-1)n+(j-1) +r}{r}
\binom{kn+(j-1)+1}{s+1-r} \prod_{i=0}^{n-1} (r+1+(j-1)+(k-1)i)
\nonumber
\end{eqnarray*}
as desired. Thus we have proved our formula for
$A^{(k)}_{s+1,kn+j}$ holds for all $n$ and for all $0 \leq j \leq k-2$.
\\

Next we verify that our formula holds
for $A^{(k)}_{s+1,kn+k-1}$.  By recursion
(\ref{kn+k-1:kn+krec}), we have
\begin{equation}\label{genform9}
(1+t)A^{(k)}_{t,kn+k-1} = A^{(k)}_{t,kn+k} - (kn+k-t)A^{(k)}_{t-1,kn+k-1}.
\end{equation}
Putting $t =s +1$ in  (\ref{genform9}), we have that
\begin{equation}\label{genform10}
(s+2)A^{(k)}_{s+1,kn+k-1} = A^{(k)}_{s+1,kn+k} -
(kn+k-(s+1))A^{(k)}_{s,kn+k-1}.
\end{equation}
Thus using the formula for $A^{(k)}_{s+1,kn+k}$ that we just proved
and our induction hypothesis, we have
\begin{eqnarray*}
&&(s+2)A^{(k)}_{s+1,kn+k-1} = \\
&&((k-1)(n+1))!\left[ \sum_{r=0}^{s+1} (-1)^{s+1-r}
\binom{(k-1)(n+1) +r}{r}
\binom{k(n+1)+1}{s+1-r} \prod_{i=0}^n (r+1 +(k-1)i)\right] \nonumber \\
&&-(kn+k-(s+1)) ((k-1)n + (k-1))! \times \nonumber \\
&& \ \ \ \ \ \ \left[ \sum_{r=0}^{s} (-1)^{s-r} \binom{(k-1)n +k-1
+r}{r} \binom{kn +k-1 +1}{s-r} \prod_{i=0}^{n-1} (r+1+(k-1)
+(k-1)i)\right].\nonumber \end{eqnarray*} Thus
\begin{eqnarray}\label{genform12}
&&\frac{A^{(k)}_{s+1,kn+k-1}}{(k-1)n +k-1)!} = \\
&&\frac{1}{s+2}\left[ \sum_{r=0}^{s+1} (-1)^{s+1-r}
\binom{((k-1)(n+1) +r}{r} \binom{k(n+1)+1}{s+1-r} \prod_{i=0}^n (r+1
+(k-1)i)\right] \nonumber
\end{eqnarray}
\begin{eqnarray*}
&&-\frac{(kn+k-(s+1))}{s+2} \times \nonumber \\
&&\left[ \sum_{r=0}^{s} (-1)^{s-r} \binom{(k-1)n +k-1 +r}{r}
\binom{kn +k}{s-r} \prod_{i=0}^{n-1} (r+1+(k-1) +(k-1)i)\right].
\nonumber
\end{eqnarray*}

Note that $\prod_{i=0}^n (r +1 +(k-1)i) = (r+1) \prod_{i=0}^{n-1}
(r +1 +(k-1)+(k-1)i)$. Thus we can divide the RHS of
(\ref{genform12}) into two terms. The first term comes from the $r
=s+1$ in the first summand on the RHS of (\ref{genform12}) and
yields
\begin{eqnarray}\label{genform13}
&&\frac{1}{s+2} \binom{((k-1)(n+1) +s+1}{s+1} (s+2)
\prod_{i=0}^{n-1} (s+1 +1 +(k-1) +(k-1)i) = \nonumber \\
&& \binom{((k-1)(n+1) +s+1}{s+1} \prod_{i=0}^{n-1} (s+1 +1 +(k-1)
+(k-1)i).
\end{eqnarray}

The remaining terms can be organized according to the factor
$\prod_{i=0}^{n-1} (r +1 +(k-1) +(k-1)i)$ to give
\begin{eqnarray}\label{genform14}
&& \sum_{r=0}^{s} (-1)^{s+1-r} \binom{(k-1)n+(k-1) +r}{r}
 \prod_{i=0}^{n-1} (r+1 + (k-1)+(k-1)i)  \frac{1}{s+2} \times \nonumber \\
&& \left[ (1+r)\binom{kn+k+1}{s+1-r} + (kn + k
-(s+1))\binom{kn+k}{s-r}\right].
\end{eqnarray}
But
\begin{eqnarray}\label{genform15}
&&\left[ (1+r)\binom{kn+k+1}{s+1-r} + (kn + k -(s+1))\binom{kn+k}{s-r}\right]
= \\
&&\frac{(kn+k)\downarrow_{s-r}}{(s+1-r)!} \left[ (1+r)(kn+k+1) +
(s+1-r)(kn+k-(s+1))\right] = \nonumber \\
&&\frac{(kn+k)\downarrow_{s-r}}{(s+1-r)!} \left[
(s+2)(kn+k) + (r+1) - (s+1-r)(s+1)\right] = \nonumber \\
&&\frac{(kn+k)\downarrow_{s-r}}{(s+1-r)!} \left[
(s+2)(kn+k) -(s+2)(s-r)\right] = \nonumber \\
&&(s+2) \frac{(kn+k)\downarrow_{s+1-r}}{(s+1-r)!} = \nonumber \\
&&(s+2) \binom{kn+k}{s+1-r}. \nonumber
\end{eqnarray}
Substituting (\ref{genform15}) in (\ref{genform14}), we get that
(\ref{genform14}) is equivalent to
\begin{equation}\label{genform16}
\sum_{r=0}^s (-1)^{s+1-r} \binom{(k-1)n+(k-1) +r}{r}
\binom{kn+k}{s+1-r} \prod_{i=0}^{n-1} (r+1 +(k-1) + (k-1)i).
\end{equation}
Finally, combining (\ref{genform13}) and (\ref{genform16}), we obtain
\begin{eqnarray*}
&&\frac{A^{(k)}_{s+1,kn+k-1}}{((k-1)n +k-1)!} = \\
&&\sum_{r=0}^{s+1} (-1)^{s+1-r} \binom{(k-1)n+(k-1) +r}{r}
\binom{kn + (k-1) +1}{s+1-r} \prod_{i=0}^{n-1} (r+1 +(k-1) +(k-1)i)  \nonumber
\end{eqnarray*}
as desired.
\end{proof}

As a corollary to Theorem~\ref{thm:5} we get combinatorial
proofs for two special cases of the Saalsch\"utz's identity, which
in terms of generalized hypergeometric functions can be written as
$${}_3F_2\left[\begin{array}{ccc} a & b & c \\ d & e & \end{array}; 1\right]=
\frac{(d-a)_{|c|}(d-b)_{|c|}}{d_{|c|}(d-a-b)_{|c|}}$$ where
$d+e=a+b+c+1$ and $c$ is a negative integer\footnote{For the first
identity in Corollary~\ref{saalschutz0}, $a=n+1$, $b=n+1$, $c=-s$,
$d=1$, and $e=2n+2-s$; for the second identity there, $a=n+2$,
$b=n+2$, $c=-s$, $d=2$, and $e=2n+3-s$} (see~\cite{A=B} pages 43 and
126).

\begin{corollary}\label{saalschutz0} The following identities hold:
$$\binom{n}{s}^2=\sum_{r=0}^{s}(-1)^{s-r}\binom{n+r}{r}^2\binom{2n+1}{s-r};$$
$$\binom{n}{s}\binom{n+1}{s+1}=\sum_{r=0}^{s}(-1)^{s-r}\binom{n+r+1}{r}\binom{n+r+1}{r+1}\binom{2n+2}{s-r}.$$
\end{corollary}

\begin{proof} The RHS of the first identity is $A^{(2)}_{s,2n}/(n!)^2$ (we use Theorem~\ref{thm:5}
for $k=2$ and $j=0$). However, as stated in Theorem~\ref{thm:KR1} in
the introduction, we proved $A^{(2)}_{s,2n} = R_{s,2n} = (n!)^2
\binom{n}{s}^2$ in \cite{kitrem}. The RHS of the second identity is
$A^{(2)}_{s,2n+1}/(n! (n+1)!)$ and by Theorem~\ref{thm:KR1},
$A^{(2)}_{s,2n+1} = R_{s,2n+1} = n!(n+1)! \binom{n}{s}
\binom{n+1}{s+1}$.
\end{proof}

Remarkably, there is a second inclusion-exclusion type formula for
the coefficients $A^{(k)}_{s,kn+j}$ which can be obtained by
interating the recursions \ref{kn+j-1:kn+jrec} and
\ref{kn+k-1:kn+krec} starting with our formulas for
$A^{(k)}_{n,kn+j}$.

\begin{theorem}\label{thm:6}
For all $0 \leq j \leq k-1$ and $0 \leq s \leq n$,
\begin{eqnarray}\label{eq:n-s,kn+j}
&&A^{(k)}_{n-s,kn+j} = \\
&&((k-1)n+j)!\left[\sum_{r=0}^{s} (-1)^{s-r} \binom{(k-1)n+j +r}{r}
\binom{kn+j+1}{s-r} \prod_{i=1}^n (r+(k-1)i)\right]. \nonumber
\end{eqnarray}

\end{theorem}
\begin{proof}

Again we proceed by induction on $s$. By Theorem~\ref{thm:2}, we
have proved our formula for $A^{(k)}_{n-s,kn+j}$
in the case where $s = 0$  for
all $n \geq 0$ and $0 \leq j \leq k-1$.

Now assume that $s > 0$ and that the theorem hold for all $s' < s$ by
induction. Note that for $n = 0$ and $j = 0, \ldots, k-1$, our formula
asserts that
\begin{eqnarray*}
A^{(k)}_{-s,j} &=& j! \sum_{r=0}^{s} (-1)^{s-r}
\binom{j +r}{r} \binom{j+1}{s-r}\\
&=& j! \left( \sum_{n \geq 0} \binom{j+n}{n} x^n \right) \left(\sum_{ m\geq 0}
(-1)^m \binom{j+1}{m}x^m\right)|_{x^s}\\
&=& j! \frac{1}{(1-x)^{j+1}} (1-x)^{j+1}|_{x^s}=0
\end{eqnarray*}
so that our formula holds for $n =0$ for $j =0, \ldots, k-1$.

Next, by induction, assume that the our formula holds for
$s$ for $n' < n$ and $j= 0, \ldots, k-1$.
Recall by (\ref{kn+k-1:kn+krec})
\begin{equation}\label{1:n-s}
A^{(k)}_{s,kn+k} =  (1+s)A^{(k)}_{s,kn+k-1} +(kn+k-s)A^{(k)}_{s-1,kn+k-1}.
\end{equation}
Replacing $n$ by $n-1$ and $s$ by $n-s$ in (\ref{1:n-s}), we
obtain that $$A^{(k)}_{n-s,kn} = (1+n-s)A^{(k)}_{n-1
-(s-1),k(n-1)+k-1} + ((k-1)n+s)A^{(k)}_{n-1 -s,k(n-1)+k-1}.$$ Thus
by induction, we get that
\begin{eqnarray}\label{3:n-s}
&&A^{(k)}_{n-s,kn} = \nonumber \\
&&(1+n-s) ((k-1)n)! \times \nonumber \\
&&\sum_{r=0}^{s-1} (-1)^{s-1-r} \binom{(k-1)(n-1)+k-1+r}{r} \binom{k(n-1)+k-1 +1}{s-1-r} \prod_{i=1}^{n-1} (r+(k-1)i) \nonumber \\
&& - ((k-1)n+s) ((k-1)n)!\times \nonumber \\
&&\sum_{r=0}^s (-1)^{s-r} \binom{(k-1)(n-1)+k-1+r}{r} \binom{k(n-1)+k-1 +1}{s-r} \prod_{i=1}^{n-1} (r+(k-1)i). \nonumber
\end{eqnarray}
Hence
\begin{eqnarray}\label{4:n-s}
&&\frac{A^{(k)}_{n-s,kn}}{((k-1)n)!} = \\
&&(1+n-s) \sum_{r=0}^{s-1} (-1)^{s-1-r} \binom{(k-1)n+r}{r} \binom{kn}{s-1-r} \prod_{i=1}^{n-1} (r+(k-1)i) \nonumber \\
&& + ((k-1)n+s) \sum_{r=0}^s (-1)^{s-r} \binom{(k-1)n +r}{r} \binom{kn}{s-r} \prod_{i=1}^{n-1} (r+(k-1)i). \nonumber
\end{eqnarray}
Now we can divide LHS of (\ref{4:n-s}) into two parts. First $r=s$ term
in the second summand is
\begin{equation}\label{5:n-s}
((k-1)n+s) \binom{(k-1)n+s}{s} \prod_{i=1}^{n-1} (s+(k-1)i) =
\binom{(k-1)n+s}{s} \prod_{i=1}^{n} (s+(k-1)i).
\end{equation}
Then we can combine the remaining terms on the LHS of (\ref{4:n-s}) to
obtain
\begin{eqnarray}\label{6:n-s}
&&\sum_{r=0}^{s-1} (-1)^{s-r} \binom{(k-1)n+r}{r} \prod_{i=1}^n (r+(k-1)i) \times \nonumber \\
&& \ \ \frac{1}{r+(k-1)n}\left[ ((k-1)n+s)\binom{kn}{s-r} -(1+n-s) \binom{kn}{s-1-r}\right].
\end{eqnarray}
The term in square brackets in (\ref{6:n-s}) is equal to
\begin{eqnarray}\label{7:n-s}
&&\left[ ((k-1)n+s)\frac{(kn)\downarrow_{s-r}}{(s-r)!}\frac{kn+1}{kn+1}
\right. \nonumber \\
&&\left.  -(1+n-s) \frac{(kn)\downarrow_{s-1-r}}{(s-1-r)!} \frac{kn+1}{kn+1}
\frac{s-r}{s-r}\right] = \nonumber \\
&&\frac{1}{kn+1} \frac{(kn+1)\downarrow_{s-r}}{(s-r)!}
[((k-1)n+s)(kn-(s-r)+1) - (1+n-s)(s-r)] = \nonumber \\
&&\frac{1}{kn+1}\binom{kn+1}{s-r} [(kn+1)(r+(k-1)n)] = \nonumber \\
&&(r+(k-1)n)\binom{kn+1}{s-r}.
\end{eqnarray}
Thus plugging in (\ref{7:n-s}) into (\ref{6:n-s}), we see that
(\ref{6:n-s}) is equal to
\begin{equation}\label{8:n-s}
\sum_{r=0}^{s-1} (-1)^{s-r} \binom{(k-1)n+r}{r} \binom{kn+1}{s-r}\prod_{i=1}^n (r+(k-1)i).
\end{equation}
Thus we can combine (\ref{5:n-s}) into (\ref{8:n-s}) to obtain that
\begin{equation}
\frac{A^{(k)}_{n-s,kn}}{((k-1)n)!} = \sum_{r=0}^{s} (-1)^{s-r} \binom{(k-1)n+r}{r} \binom{kn+1}{s-r}\prod_{i=1}^n (r+(k-1)i).
\end{equation}
as desired.

By (\ref{kn+j-1:kn+jrec}), we have that
\begin{equation}\label{9:n-s:1}
A^{(k)}_{s,kn+j} =  (kn+j-s)A^{(k)}_{s,kn+j-1} +(s+1)A^{(k)}_{s+1,kn+j-1}
\end{equation}
for $1 \leq j \leq k-1$. Replacing $s$ by $n-s$ in
(\ref{9:n-s:1}), we obtain that $$A^{(k)}_{n-s,kn+j} =
((k-1)n+j+s)A^{(k)}_{n-s,kn+j-1}
+(n-s+1)A^{(k)}_{n-(s-1),kn+j-1}.$$ We can assume by induction
that our formula holds for $A^{(k)}_{n-s,kn+j-1}$ so that
\begin{eqnarray}\label{10:n-s}
&& A^{(k)}_{n-s,kn+j} =  \nonumber \\
&& ((k-1)n+j+s)((k-1)n+j-1)! \times \nonumber \\
&&\sum_{r=0}^s (-1)^{s-r}\binom{(k-1)n+j-1+r}{r}
\binom{kn+j-1+1}{s-r} \prod_{i=1}^n (r + (k-1)i) \nonumber \\
&& +(n-s+1)((k-1)n+j-1)! \times \nonumber \\
&&\sum_{r=0}^{s-1} (-1)^{s-1-r}\binom{(k-1)n+j-1+r}{r}
\binom{kn+j-1+1}{s-1-r} \prod_{i=1}^n (r + (k-1)i) \nonumber.
\end{eqnarray}
Thus
\begin{eqnarray}\label{11:n-s}
&& \frac{A^{(k)}_{n-s,kn+j}}{((k-1)n+j-1)!} =  \\
&& ((k-1)n+j+s)\sum_{r=0}^s (-1)^{s-r}\binom{(k-1)n+j-1+r}{r}
\binom{kn+j}{s-r} \prod_{i=1}^n (r + (k-1)i) \nonumber \\
&& +(n-s+1)\sum_{r=0}^{s-1} (-1)^{s-1-r}\binom{(k-1)n+j-1+r}{r}
\binom{kn+j}{s-1-r}\prod_{i=1}^n (r + (k-1)i) \nonumber.
\end{eqnarray}
Again, we can divide the LHS of (\ref{11:n-s}) into two terms.
First term coming from the $r=s$ of the first summand is
\begin{eqnarray}\label{12:n-s}
&&((k-1)n+j+s) \binom{(k-1)n+j-1+s}{s} \prod_{i=1}^n (s+(k-1)i) = \\
&& ((k-1)n+j) \binom{(k-1)n+j+s}{s}\prod_{i=1}^n (s+(k-1)i)
\nonumber.
\end{eqnarray}

Next the remaining terms on the LHS of (\ref{11:n-s}) can be
combined into
\begin{eqnarray}\label{13:n-s}
&&\sum_{r=0}^{s-1} (-1)^{s-r} \prod_{i=1}^n (s+(k-1)i) \times \nonumber \\
&&\left[ ((k-1)n+j+s)\binom{(k-1)n+j-1+r}{r}
\binom{kn+j}{s-r} - \right.\nonumber \\
&&\left. (n-s+1)\binom{(k-1)n+j-1+r}{r}
\binom{kn+j}{s-1-r}\right].
\end{eqnarray}
Now the term in the square brackets in (\ref{13:n-s}) is equal to
\begin{eqnarray}\label{14:n-s}
&&((k-1)n+j+s) \frac{((k-1)n+j+r-1)\downarrow_r}{r!}
\frac{(k-1)n+j+r}{(k-1)n+j+r}
\frac{(kn+j)\downarrow_{s-r}}{(s-r)!} \frac{kn+j+1}{kn+j+1} \nonumber \\
&& - (n-s+1) \frac{((k-1)n+j+r-1)\downarrow_r}{r!}
\frac{(k-1)n+j+r}{(k-1)n+j+r}
\frac{(kn+j)\downarrow_{s-1-r}}{(s-1-r)!}
\frac{kn+j+1}{kn+j+1} \frac{s-r}{s-r} = \nonumber \\
&& \frac{((k-1)n+j+r)\downarrow_r}{r!}
\frac{(kn+j+1)\downarrow_{s-r}}{(s-r)!}  \times \nonumber \\
&&\left[ ((k-1)n+j+s) \frac{(k-1)n+j}{(k-1)n+j+r} \frac{kn+j-(s-r)+1}{kn+j+1} \right. \nonumber \\
&&\left.- (n-s+1)\frac{(k-1)n+j}{(k-1)n+j+r} \frac{(s-r)}{kn+j+1}
\right]=
\nonumber \\
&& \binom{(k-1)n+j+r}{r} \binom{kn+j+1}{s-r}
\frac{(k-1)n+j}{((k-1)n+j+r)(kn+j+1)} \times \nonumber \\
&&\left[((k-1)n+j+s)(kn+j-(s-r)+1) - (n-s+1)(s-r)\right]= \nonumber \\
&& \binom{(k-1)n+j+r}{r} \binom{kn+j+1}{s-r}
\frac{(k-1)n+j}{((k-1)n+j+r)(kn+j+1)} \times \nonumber \\
&&\left[ (k-1)n+j+r)(kn+j+1) \right]= \nonumber \\
&&((k-1)n+j) \binom{(k-1)n+j+r}{r} \binom{kn+j+1}{s-r}.
\end{eqnarray}
Here the fact that $((k-1)n+j+s)(kn+j-(s-r)+1) - (n-s+1)(s-r) =
(k-1)n+j+r)(kn+j+1)$ can easily be verified in your favorite computer
algebra system.
Thus substituting in the result of (\ref{14:n-s}) into
(\ref{13:n-s}), we see that (\ref{13:n-s}) is equal to
\begin{equation}\label{15:n-s}
((k-1)n+j) \sum_{r=0}^{s-1} (-1)^{s-r} \binom{(k-1)n+j+r}{r}
\binom{kn+j+1}{s-r}\prod_{i=1}^n (r+(k-1)i).
\end{equation}
Combining (\ref{12:n-s}) and (\ref{15:n-s}), we see that
\begin{eqnarray}\label{16:n-s}
&&\frac{A^{(k)}_{n-s,kn+j}}{((k-1)n+j-1)!} = \\
&&((k-1)n+j) \sum_{r=0}^{s} (-1)^{s-r} \binom{(k-1)n+j+r}{r}
\binom{kn+j+1}{s-r}\prod_{i=1}^n (r+(k-1)i) \nonumber
\end{eqnarray}
or, equivalently, that
 \begin{eqnarray}\label{17:n-s}
&&A^{(k)}_{n-s,kn+j} = \nonumber \\
&&((k-1)n+j)! \sum_{r=0}^{s} (-1)^{s-r} \binom{(k-1)n+j+r}{r}
\binom{kn+j+1}{s-r}\prod_{i=1}^n (r+(k-1)i) \nonumber
\end{eqnarray}
as desired.
This completes our induction and hence we have established
our formulas for $A^{(k)}_{n-s,kn+j}$ for all $0 \leq j \leq k-1$  and $n \geq 0$.
\end{proof}

Note if we compare the formulas for $A^{(k)}_{n-s,kn+j}$ from
Theorem \ref{thm:5} and from Theorem \ref{thm:6}, we obtain the
following identities
\begin{corollary}\label{cor1} For all $0 \leq j \leq k-1$ and
$0 \leq s \leq n$,
\begin{eqnarray}\label{s=n-s}
&&\sum_{r=0}^{s} (-1)^{s-r}
\binom{(k-1)n+j+r}{r} \binom{kn+j+1}{s-r}\prod_{i=1}^n (r+(k-1)i) = \nonumber \\
&&\sum_{r=0}^{n-s} (-1)^{n-s-r} \binom{(k-1)n+j+r}{r}
\binom{kn+j+1}{n-s-r}\prod_{i=0}^{n-1} (r+1+j+(k-1)i). \nonumber
\end{eqnarray}
\end{corollary}

For example, the $s =0$ of Corollary \ref{cor1}, gives the following
indentities.  For all $n$ and for all $0 \leq j \leq k-1$,
\begin{equation*}
(k-1)^n (n!) =
\sum_{r=0}^n (-1)^{n-r} \binom{(k-1)n + j +r}{r} \binom{kn+j+1}{n-r}
\prod_{i=0}^{n-1} (r+1+j + (k-1)i).
\end{equation*}

Similarly, the $s =1$ of Corollary \ref{cor1}, gives the following
indentities. For all $n$ and for all $0 \leq j \leq k-1$,
\begin{eqnarray*}
&&((k-1)n+j+1)\prod_{i=0}^n (1+(k-1)i) - (kn+j+1) (k-1)^n (n!) = \\
&&\sum_{r=0}^{n-1} (-1)^{n-1-r} \binom{(k-1)n + j +r}{r} \binom{kn+j+1}{n-1-r}
\prod_{i=0}^{n-1} (r+ n + (k-1)i). \nonumber
\end{eqnarray*}

\section{Properties of $B_n^{(k)}(x,z)$}\label{section3}

For $0 \leq j \leq k-2$, let $\Theta_{kn+j}$ be the operator that
sends $z^0x^s$ to $(1 + s+ (k-1)n +j)z^0 x^s + (n-s)z^0 x^{s+1}$
and $z^1x^s$ to $(1 + s+ (k-1)n +j)z^1 x^s + (n-s-1)z^1 x^{s+1} +
z^0 x^{s+1}$. Also let $\Psi_{kn+k-1}$ be the operator that sends
$z^0x^s$ to $(s+ (k-1)(n +1))z^0 x^s + z^1 x^s+ (n-s)z^0 x^{s+1}$
and $z^1x^s$ to $(1 + s+ (k-1)(n +1))z^1 x^s + (n-s)z^1 x^{s+1}$.
Then we have the following.

\begin{theorem}\label{thm:B1}
For any $k \geq 2$ and $n \geq 0$,
\begin{enumerate}
\item $B^{(k)}_1(x,z) = 1$,

\item $B^{(k)}_{kn+j+1}(x,z) = \Theta_{kn+j}(B^{(k)}_{kn+j}(x,z))$ for $0 \leq j \leq k-2$, and

\item  $B^{(k)}_{kn+k}(x,z) =
\Psi_{kn+k-1}(B^{(k)}_{kn+k-1}(x,z))$.
\end{enumerate}
\end{theorem}
\begin{proof}
Part (1) is easy to see.

For part (2), suppose we are given a permutation $\sg = \sg_1
\cdots \sg_{kn+j}$ such that $\overrightarrow{des}_{kN}(\sg) =s $
and $\sg_1 \notin kN$ so that such a $\sg$ gives rise to a factor
of $z^0 x^s$ in $B^{(k)}_{kn+j}(x,z)$. Then for each $i$ such that
$\sg_i > \sg_{i+1}$ and $\sg_{i+1} \in kN$, inserting $kn+j+1$ in
position $i$ will result in a permutation $\sg^{(i)}$ such that
$\overrightarrow{des}_{kN}(\sg^{(i)}) =s$ and  $\sg^{(i)}_1 \notin
kN$. Similarly if $i = kn+j$ or $i$ is such that $\sg_{i+1} \notin
kN$, then $\overrightarrow{des}_{kN}(\sg^{(i)}) =s$ and
$\sg^{(i)}_1 \notin kN$. Finally if $i$ is such that $i \notin
\overrightarrow{Des}_{kN}(\sg)$ and $\sg_{i+1} \in kN$, then
$\overrightarrow{des}_{kN}(\sg^{(i)}) =s+1$ and $\sg^{(i)}_1
\notin kN$. Thus the collection of $\{\sg^{(i)}: i =0, \ldots,
kn+j\}$ contributes a factor of  $(1 + s+ (k-1)n +j)z^0 x^s +
(n-s)z^0 x^{s+1}$ to $B^{(k)}_{kn+j+1}(x,z)$.

Next suppose we are given a permutation $\sg = \sg_1 \cdots
\sg_{kn+j}$ such that $\overrightarrow{des}_{kN}(\sg) =s $ and
$\sg_1 \in kN$ so that such a $\sg$ gives rise to a factor of $z^1
x^s$ in $B^{(k)}_{kn+j}(x,z)$. Then for each $i$ such that $\sg_i
> \sg_{i+1}$ and $\sg_{i+1} \in kN$, inserting $kn+j+1$ in
position $i$ will result in a permutation $\sg^{(i)}$ such that
$\overrightarrow{des}_{kN}(\sg^{(i)}) =s$ and  $\sg^{(i)}_1 \in
kN$. Similarly if $i = kn+j$ or $i$ is such that $\sg_{i+1} \notin
kN$, then $\overrightarrow{des}_{kN}(\sg^{(i)}) =s$ and
$\sg^{(i)}_1 \in kN$. If $i > 0$ is such that $i \notin
\overrightarrow{Des}_{kN}(\sg)$ and $\sg_{i+1} \in kN$, then
$\overrightarrow{des}_{kN}(\sg^{(i)}) =s+1$ and $\sg^{(i)}_1 \in
kN$. Finally if $i =0$, then, since $\sg_1 \in kN$,
$\overrightarrow{des}_{kN}(\sg^{(0)}) =s+1$ and $\sg^{(0)}_1
\notin kN$. Thus the collection of $\{\sg^{(i)}: i =0, \ldots,
kn+j\}$ contributes a factor of  $(1 + s+ (k-1)n +j)z^1 x^s +
(n-s-1)z^1 x^{s+1} +z^0 x^{s+1}$ to $B^{(k)}_{kn+j+1}(x,z)$ in
this case.

For part (3), suppose $\sg = \sg_1
\cdots \sg_{kn+k-1}$ is a permutation
such that $\overrightarrow{des}_{kN}(\sg) =s
$ and $\sg_1 \notin kN$ so that such a $\sg$ gives rise to a
factor of $z^0 x^s$ in $B^{(k)}_{kn+k-1}(x,z)$. Then for each $i$
such that $\sg_i > \sg_{i+1}$ and $\sg_{i+1} \in kN$, inserting
$kn+k$ in position $i$ will result in a permutation $\sg^{(i)}$
such that $\overrightarrow{des}_{kN}(\sg^{(i)}) =s$ and
$\sg^{(i)}_1 \notin kN$. Similarly if $i = kn+k-1$ or $i > 0$ is
such that $\sg_{i+1} \notin kN$, then
$\overrightarrow{des}_{kN}(\sg^{(i)}) =s$ and $\sg^{(i)}_1 \notin
kN$. Since $\sg_1 \notin kN$, inserting $kn+k$ in position 0 will
result in a permutation such that
$\overrightarrow{des}_{kN}(\sg^{(0)}) =s$ and $\sg^{(i)}_1 \in
kN$. Finally if $i$ is such that $i \notin
\overrightarrow{Des}_{kN}(\sg)$ and $\sg_{i+1} \in kN$, then
$\overrightarrow{des}_{kN}(\sg^{(i)}) =s+1$ and $\sg^{(i)}_1
\notin kN$. Thus the collection of $\{\sg^{(i)}: i =0, \ldots,
kn+k-1\}$ contributes a factor of  $(s+ (k-1)n +k-1)z^0 x^s + z^1
x^s + (n-s)z^0 x^{s+1}$ to $B^{(k)}_{kn+k}(x,z)$.

Next suppose $\sg = \sg_1 \cdots
\sg_{kn+k-1}$ is a permutation
such that $\overrightarrow{des}_{kN}(\sg) =s $ and
$\sg_1 \in kN$ so that such a $\sg$ gives rise to a factor of $z^1
x^s$ in $B^{(k)}_{kn+k-1}(x,z)$. Then for each $i$ such that
$\sg_i > \sg_{i+1}$ and $\sg_{i+1} \in kN$, inserting $kn+k$ in
position $i$ will result in a permutation $\sg^{(i)}$ such that
$\overrightarrow{des}_{kN}(\sg^{(i)}) =s$ and $\sg^{(i)}_1 \in
kN$. Similarly if $i = kn+k-1$ or $i$ is such that $\sg_{i+1}
\notin kN$, then $\overrightarrow{des}_{kN}(\sg^{(i)}) =s$ and
$\sg^{(i)}_1 \in kN$. Finally if $i$ is such that $\sg_{i+1} \in
kN$, then $\overrightarrow{des}_{kN}(\sg^{(i)}) =s+1$ and
$\sg^{(i)}_1 \in kN$. Thus the collection of $\{\sg^{(i)}: i =0,
\ldots, kn+k-1\}$ contributes a factor of  $(s+ (k-1)n +k)z^1 x^s
+ (n-s)z^1 x^{s+1}$ to $B^{(k)}_{kn+k}(x,z)$ in this case.
\end{proof}

One can also express the actions $\Theta_{kn+j}$ and $\Psi_{kn+k-1}$ in
terms of partial differential operators. That is, it is easy to verify
that we have the following.

\begin{theorem}\label{thm:B2}
For any $k \geq 2$ and $n \geq 0$,
\begin{enumerate}
\item $B^{(k)}_1(x,z) = 1$,

\item $B^{(k)}_{kn+j+1}(x,z) =
\left( x(1-x)\frac{\partial}{\partial x} + x(1-z)\frac{\partial}{\partial z}
+ nx + (1 +(k-1)n +j) \right)(B^{(k)}_{kn+j}(x,z))$ for $0 \leq j \leq k-2$, and

\item  $B^{(k)}_{kn+k}(x,z) = \left(
x(1-x)\frac{\partial}{\partial x} + z(1-z)\frac{\partial}{\partial
z} + nx + z +(k-1)(n +1) \right)(B^{(k)}_{kn+k-1}(x,z))$.
\end{enumerate}
\end{theorem}
One can use Theorem \ref{thm:B2} to generate expressions for
$B^{(k)}_n(x,z)$ for small values of $k$ and $n$. Here are some
initial values for
$B^{(3)}_n(x,z)$.\\
\begin{description}

\item $B^{(3)}_1(x,z) =1$.

\item $B^{(3)}_2(x,z) =2$.

\item $B^{(3)}_3(x,z) =4+2z$.

\item $B^{(3)}_4(x,z) =12+6z+ 6x$.

\item $B^{(3)}_5(x,z) =48+24z+48x$.

\item $B^{(3)}_6(x,z) =192 + 168z+288x + 72xz$.

\item $B^{(3)}_7(x,z) =960+840z+ 2280x +600xz + 360x^2$.

\item $B^{(3)}_8(x,z) =5760+5040z + 18720x + 5040z + 5760x^2$.

\item $B^{(3)}_9(x,z) = 34560+41040z+142560x + 69120xz + 64800 x^2 +
10800x^2z$.

\item $B^{(3)}_{10}(x,z) = 241920 + 287280z + 1285200x + 635040xz +
937440x^2 + 166320 x^2z + 75600 x^3$.

\item $B^{(3)}_{11}(x,z) = 1935360 + 2298240 z+ 12579840 x + 6289920xz +
12579840 x^2 + 2298240 x^2 z + 1935360x^3$.

\item $B^{(3)}_{12}(x,z) = 15482880 + 22619520z + 119024640x + 82373760 xz +
150958080x^2 + 50440320 x^2z + 33868800x^3 + 4233600x^3z$.

\end{description}

Recall that for all $n$, we let $B^{(k)}_n(x,z) = \sum_{\sg \in
\S_n} x^{\overrightarrow{des}_{kN}(\sg)}z^{\chi(\sg_1 \in kN)} =
\sum_{j=0}^{\lfloor \frac{n}{k} \rfloor} \sum_{i=0}^1
B^{(k)}_{i,j,n} z^ix^j$ and $B^{(k)}_{s,kn+j} = B^{(k)}_{0,s,kn+j}
+ B^{(k)}_{1,s,kn+j}$, for all $n$ and $j$ where $0 \leq j \leq
k-1$. It is easy to see that Theorem~\ref{thm:B1} implies to the
following recursions.

For $0 \leq j \leq k-2$ and all $n \geq 0$,
\begin{equation}\label{rec:B0j}
B^{(k)}_{0,s,kn+j+1} = (1+s + (k-1)n +j)B^{(k)}_{0,s,kn+j} + (n-s+1)B^{(k)}_{0,s-1,kn+j} +
B^{(k)}_{1,s-1,kn+j},
\end{equation}
\begin{equation}\label{rec:B1j}
B^{(k)}_{1,s,kn+j+1} = (1+s + (k-1)n +j)B^{(k)}_{1,s,kn+j} + (n-s)B^{(k)}_{1,s-1,kn+j},
\end{equation}
and
\begin{equation}\label{rec:Cj}
B^{(k)}_{s,kn+j+1} = (1+s + (k-1)n +j)B^{(k)}_{s,kn+j} +
(n-s+1)B^{(k)}_{s-1,kn+j}.
\end{equation}
Similarly for all $n \geq 0$,
\begin{equation}\label{rec:B0k-1}
B^{(k)}_{0,s,kn+k} = (s + (k-1)(n +1))B^{(k)}_{0,s,kn+k-1} +
(n-s+1)B^{(k)}_{0,s-1,kn+k-1},
\end{equation}
\begin{equation}\label{rec:B1k-1}
B^{(k)}_{1,s,kn+k} = (1+s + (k-1)(n +1))B^{(k)}_{1,s,kn+k-1} +
(n-s+1)B^{(k)}_{1,s-1,kn+k-1} + B^{(k)}_{0,s,kn+k-1},
\end{equation}
and
\begin{equation}\label{rec:Ck-1}
B^{(k)}_{s,kn+k} = (1+s + (k-1)(n +1))B^{(k)}_{s,kn+k-1} +
(n-s+1)B^{(k)}_{s-1,kn+k-1}.
\end{equation}
Note that~(\ref{rec:Ck-1}) can be obtained from~(\ref{rec:Cj}) by
plugging in $j=k-1$. Next we shall consider some special values of
$B^{(k)}_{i,j,n}$ and $B^{(k)}_{j,n}$.

\begin{theorem}\label{thm:B3}
For all $n \geq 0$, $k \geq 2$, and $0 \leq j \leq k-1$,
\begin{enumerate}
\item $B^{(k)}_{0,kn+j} = ((k-1)n+j)! \prod_{i=1}^n (1 + (k-1)i)$,

\item $B^{(k)}_{0,0,kn+j} = ((k-1)n+j)! \prod_{i=1}^n ((k-1)i) =
((k-1)n+j)! (k-1)^n (n!)$, and

\item $B^{(k)}_{1,0,kn+j} = ((k-1)n+j)! \left( (\prod_{i=1}^n (1+ (k-1)i)) -
(k-1)^n (n!) \right)$.
\end{enumerate}
\end{theorem}
\begin{proof}
It is easy to see that (1), (2) and (3) hold for $n=1$.

To prove (1), fix $j$ such that $0 \leq j \leq k-1$ and suppose
that $\sg = \sg_1 \cdots \sg_{kn+j}$ is such that
$\overrightarrow{des}_{kN}(\sg) =0$. Then we can factor any such
permutation into blocks by reading the permutation from right to
left  and cutting before each  number which is not equal to $0 \mod
k$. For example if $j=0$, $k=3$, and $\sg =
11~1~3~6~2~4~5~7~8~9~10~12$, then the blocks of $\sg$ would be
$11$, $1~3~6~$, $2$,  $4$, $5$, $7$,  $8~9$, $10~12$.

There may be a block of numbers which are equivalent to $0 \mod k$
at the start which are arranged in increasing order. We call this
initial block the $0$-th block. Every other block must start with
a number $sk+i$ where $0 \leq s \leq n$ and $1 \leq i \leq k-1$
which can be followed  by any subset of numbers which are
equivalent to $0 \mod k$ and which are greater than $sk+i$
arranged in increasing order. We call such a block the $(sk+i)$-th
block. It is then easy to see that there are $\prod_{i=1}^{n}
(1+i(k-1))$ ways to put the numbers $k, 2k, \ldots, nk$ into the
blocks.  That is, $k$ must be placed  in one of the  blocks $0,1,
\ldots, k-1$ so that there are $1+(k-1)$ choices for the block to
place $k$.  Then $2k$ must be placed in one of the blocks $0,1,
\ldots, k-1, k+1, \ldots, k+k-1$ so that there are $1 +2 (k-1)$
choices for the block that contains $2k$. More generally, $ik$
must be placed in one of the blocks $0,1, \ldots, k-1, k+1,
\ldots, k+k-1, \ldots, (i-1)k+1, \ldots, (i-1)k +k-1$ so that
there are $1 +i (k-1)$ choices for the block that contains $ik$.
Once we have arranged the numbers which are equivalent to $0 \mod
k$ into blocks, it is easy to see that we can arrange blocks
$sk+i$ where $0 \leq s \leq n-1$ and $1 \leq i \leq k-1$ plus the
blocks $kn+1, \ldots, kn+j$ in any order and still get a
permutation $\sg$ with $\overrightarrow{des}_{kN}(\sg)= 0$. It
thus follows that there are $((k-1)n+j)! \prod_{i=1}^{n} (1 +
i(k-1))$ such permutations so that $B^{(k)}_{0,kn+j} = ((k-1)n+j)!
\prod_{i=1}^{n} (1 + i(k-1))$.

To prove our formula for $B^{(k)}_{0,0,kn+j}$, we can essentially
use the same argument except that in this case we cannot put any
elements in the $0$-th block since we want to produce permutations
that do not start with an element which is $\equiv 0 \mod k$. Thus
there are only $\prod_{i=1}^n (k-1)i = (k-1)^n n!$ ways to put the
elements which are $\equiv 0 \mod k$ in this case. Once again,
after we have arranged the numbers which are equivalent to $0 \mod
k$ into blocks, it is easy to see that we can arrange blocks
$sk+i$ where $0 \leq s \leq n-1$ and $1 \leq i \leq k-1$ plus the
blocks $kn+1, \ldots, kn+j$ in any order and still get a
permutation $\sg$ with $\overrightarrow{des}_{kN}(\sg)= 0$ and
$\sg_1 \notin kN$. It thus follows that there are $((k-1)n+j)!
\prod_{i=1}^{n} (i(k-1))$ such permutations so that
$B^{(k)}_{0,0,kn+j}=  ((k-1)n+j)! \prod_{i=1}^{n} (i(k-1))$.

Clearly, part (3) immediately follows from parts (1) and (2).

\end{proof}

\begin{theorem}\label{thm:B4}
We have $B^{(k)}_{n,kn} = B^{(k)}_{0,n,kn} = 0$, and for all $n
\geq 0$, $k \geq 2$, and $1 \leq j \leq k-1$,
\begin{enumerate}
\item $B^{(k)}_{n,kn+j} = ((k-1)n+j)! \prod_{i=0}^{n-1} (j +
(k-1)i)$,

\item $B^{(k)}_{0,n,kn+j} = ((k-1)n+j)!\prod_{i=0}^{n-1} (j +
(k-1)i)$, and

\item $B^{(k)}_{1,n,kn+j} = 0$.
\end{enumerate}
\end{theorem}
\begin{proof}

It is easy to see that if a permutation $\sg$ starts with an
element which is $\equiv 0 \mod k$, then it cannot be the case
that $\overrightarrow{des}_{kN}(\sg) =n$. Thus $B^{(k)}_{1,n,kn+j}
=0$ and $B^{(k)}_{n,kn+j} = B^{(k)}_{0,n,kn+j}$. Moreover,
$B^{(k)}_{0,n,kn} = 0$ since for permutations $\sg \in \S_{kn}$,
$kn$ cannot be the second element of a descent because it is the
largest element in the permutation. So we need only prove~(1).

To prove (1), fix $j$ such that $1 \leq j \leq k-1$ and suppose
that $\sg = \sg_1 \cdots \sg_{kn+j}$ is such that
$\overrightarrow{des}_{kN}(\sg) =n$. Then, as above, we can factor
any such permutation into blocks by reading the permutation from
right to left and cutting before each number which is not equal to
0 mod $k$. One can see that, unlike the case where
$\overrightarrow{des}_{kN}(\sg) =0$, there can be no numbers which
are equivalent to 0 mod $k$ at the start so that the $0$-th block
must be empty.  Next if $0 \leq s \leq n-1$ and $1 \leq i \leq
k-1$ and there are numbers which are equivalent to $0 \mod k$ in
the $(sk+i)$-th block, then those numbers must all be less  than
$sk+i$ and they must be arranged in decreasing order. Finally, if
there are numbers in any of the blocks $kn+1, \ldots, kn+j$
contains numbers which are equivalent to $0 \mod k$, then those
numbers must be arranged in decreasing order. It is then easy to
see that there are $\prod_{i=0}^{n-1}(j+(k-1)i)$ ways to put the
numbers $k, 2k, \ldots, nk$ into blocks.  That is, $kn$ can only
be placed in one of the blocks $kn+1, \ldots, kn+j$ so that there
are only $j$ possible choices for the blocks containing $kn$. Then
$k(n-1)$ can only be placed in one of the blocks $k(n-1)+1,
\ldots, k(n-1) +k-1, kn+1, \ldots, kn+j$  so that there are only
$j +(k-1)$ possible choices for the blocks containing $k(n-1)$. In
general, $k(n-i)$ can only be placed in one of the blocks
$k(n-i)+1, \ldots, k(n-i) +k-1, \ldots, k(n-1)+1, \ldots, k(n-1)
+k-1, kn+1, \ldots, kn+j$ so that there are only $j +i(k-1)$
possible choices for the blocks containing $k(n-i)$.
 Once we have arranged the numbers which are equivalent to
$0 \mod k$ into blocks, it is easy to see that we can arrange
blocks $sk+i$ where $0 \leq s \leq n-1$ and $1 \leq i \leq k-1$
plus the blocks $kn+1, \ldots, kn+j$ in any order an still get a
permutation $\sg$ with $\overrightarrow{des}_{kN}(\sg)= n$. It
thus follows that there are $(n(k-1)+j)! \prod_{i=0}^{n-1}
j+(k-1)i$ such permutations.
\end{proof}

\begin{theorem}\label{thm:B4:1}
For all $n \geq 0$, $k \geq 2$, and $0 \leq j \leq k-1$,
\begin{equation}\label{B1n-1}
B^{(k)}_{1,n-1,kn+j} = ((k-1)n+j)! \sum_{p=0}^{n-1} \left(
\prod_{i=0}^{p-1} (j+(k-1)i) \right) \left( \prod_{i=p+1}^{n-1}
(1+j+(k-1)i) \right).
\end{equation}
\end{theorem}
\begin{proof}

Fix $j$ such that $0 \leq j \leq k-1$ and suppose that $\sg =
\sg_1 \cdots \sg_{kn+j}$ is such that
$\overrightarrow{des}_{kN}(\sg) =n-1$. Since we are computing
$B^{(k)}_{1,n-1,kn+j}$, we need only consider permutations  $\sg$
such that  $\sg_1 \in kN$. Now if $\overrightarrow{des}_{kN}(\sg)
=n-1$, then $\sg_1$ is the only element which is equivalent to $0
\mod k$ which is not the second element of a descent pair. Thus we
will classify our permutations by the first element. That is,
suppose $\sg_1 = (n-p)k$ for some $0 \leq p \leq n-1$. Then, as
above, we can factor any such permutation into blocks by reading
the permutation from right to left and cutting before each number
which is not equal to 0 mod $k$. Thus the  $0$-th block must start
with $(n-p)k$ which could be followed by elements which are
equivalent to $0 \mod k$ and which are less than $(n-p)k$ arranged
in decreasing order.  Next if $0 \leq s \leq n$ and $1 \leq i \leq
k-1$ and there are numbers which are equivalent to $0 \mod k$ in
the $(sk+i)$-th block, then those numbers must all be less than
$sk+i$ and they must be arranged in decreasing order. Finally, if
there are numbers in any of the blocks $kn+1, \ldots, kn+j$ which
are equivalent to $0 \mod k$, then those numbers must be arranged
in decreasing order.  We claim that there are
$\left(\prod_{i=0}^{p-1} (j+(k-1)i)\right) \left(
\prod_{i=p+1}^{n-1} (1+j+(k-1)i)\right)$ ways to put the numbers
$k, 2k, \ldots, (n-p-1)k,(n-p+1)k,\ldots,nk$ into blocks.  That
is, suppose that $p \geq 1$ and consider the possibilities for the
placements of $kn, k(n-1), \ldots, k(n-p+1)$ into blocks. Clearly,
$kn$ can only be placed in one of the blocks $kn+1, \ldots, kn+j$
so that there are only $j$ possible choices for the blocks
containing $kn$. Then $k(n-1)$ can only be placed in one of the
blocks $k(n-1)+1, \ldots, k(n-1) +k-1, kn+1, \ldots, kn+j$  so
that there are only $j +(k-1)$ possible choices for the blocks
containing $k(n-1)$. In general, $k(n-i)$ can only be placed in
one of the blocks $k(n-i)+1, \ldots, k(n-i) +k-1, \ldots,
k(n-1)+1, \ldots, k(n-1) +k-1, kn+1, \ldots, kn+j$ so that there
are only $j +i(k-1)$ possible choices for the blocks containing
$k(n-i)$. Thus the total number of choices for the placements of
$kn, k(n-1), \ldots, k(n-p+1)$ into blocks is $\prod_{i=0}^{p-1}
(j+(k-1)i)$. Next consider the possibilities for the placements of
$k(n-p-1), k(n-p-2), \ldots, k$ into blocks. For $i \geq 1$,
$k(n-p-i)$ can be placed in the 0-th block or any of the blocks
$k(n-p-i)+1, \ldots, k(n-p-i) +k-1, \ldots, k(n-1)+1, \ldots,
k(n-1) +k-1, kn+1, \ldots, kn+j$ so that there are  $1+ j
+(p+i)(k-1)$ possible choices for the blocks containing
$k(n-p-i)$. Thus the total number of choices for the placements of
$k(n-p-1), k(n-p-2), \ldots, k$ into blocks is
$\prod_{i=p+1}^{n-1} (1+j+(k-1)i)$. A similar argument will show
that if $\sg_1 = kn$, then there $\prod_{i=1}^{n-1} (1+ j +
(k-1)i)$ ways to partition the elements which are equivalent to $0
\mod k$ into the blocks. After we have partitioned the  numbers
which are equivalent to $0 \mod k$ into their respective blocks,
we must arrange the numbers which are equivalent to $0 \mod k$ in
each block in decreasing order so that there are a total
$\left(\prod_{i=0}^{p-1} (j+(k-1)i)\right) \left(
\prod_{i=p+1}^{n-1} (1+j+(k-1)i)\right)$
 ways to
partition the numbers which are equivalent to $0 \mod k$ into the
blocks for a permutation $\sg$ such that $\sg_1 = k(n-p)$ and
$\overrightarrow{des}_{kN}(\sg)= n-1$ for $p =0, \ldots, n-1$.
Once we have arranged the numbers which are equivalent to $0 \mod
k$ into blocks, it is easy to see that we can arrange blocks
$sk+i$ where $0 \leq s \leq n-1$ and $1 \leq i \leq k-1$ plus the
blocks $kn+1, \ldots, kn+j$ in any order and still get a
permutation $\sg$ with $\overrightarrow{des}_{kN}(\sg)= n-1$. It
thus follows that there are $(n(k-1)+j)! \left(\prod_{i=0}^{p-1}
(j+(k-1)i)\right) \left( \prod_{i=p+1}^{n-1} (1+j+(k-1)i)\right)$
such permutations that start with $k(n-p)$, so
$B^{(k)}_{1,n-1,kn+j} = ((k-1)n+j)! \sum_{p=0}^{n-1} \left(
\prod_{i=0}^{p-1} (j+(k-1)i) \right) \left( \prod_{i=p+1}^{n-1}
(1+j+(k-1)i) \right)$.
\end{proof}

We let $\Omega(k,n,j) = \sum_{p=0}^{n-1} \left( \prod_{i=0}^{p-1}
(j+(k-1)i) \right) \left( \prod_{i=p+1}^{n-1} (1+j+(k-1)i)
\right)$. It turns out that there is another expression for
$\Omega(k,n,j)$ which will be useful for our later developments.

\begin{lemma}\label{lem:omega}
For all $n \geq 1$, $k \geq 2$, and $r \geq 0$,
\begin{eqnarray}\label{omega1}
\Omega(k,n,r) &=& \sum_{p=0}^{n-1} \left( \prod_{i=0}^{p-1}
(r+(k-1)i) \right)
\left( \prod_{i=p+1}^{n-1} (1+r+(k-1)i) \right) \nonumber \\
 &=& \left( \prod_{i=0}^{n-1} (1+r+(k-1)i) \right) -
\left(\prod_{i=0}^{n-1} (r+(k-1)i) \right).
\end{eqnarray}
\end{lemma}
\begin{proof}

Note that for $n =1$, $\sum_{p=0}^{1-1} \left( \prod_{i=0}^{p-1}
(r+(k-1)i) \right) \left( \prod_{i=p+1}^{1-1} (1+r+(k-1)i) \right)
=1$ since both products involved in the sum are empty. However,
$\left( \prod_{i=0}^{1-1} (1+r+(k-1)i) \right) -
\left(\prod_{i=0}^{1-1} (r+(k-1)i) \right) = (1+r) -r =1$ so that
our formula holds for all $r$ when $n=1$.

Next suppose that $m >1$ and that our formula holds for all $r$ when
$n = m-1$. Then

\begin{eqnarray*}
&&\left( \prod_{i=0}^{m-1} (1+r+(k-1)i) \right) - \left(\prod_{i=0}^{m-1} (r+(k-1)i) \right) = \\
&&(1+r)\left( \prod_{i=1}^{m-1} (1+r+(k-1)i) \right) - r
\left(\prod_{i=1}^{m-1} (r+(k-1)i) \right) =
\end{eqnarray*}
\begin{eqnarray*}
&&\left(\prod_{i=1}^{m-1} (1+r+(k-1)i) \right) +\\
&&r \left( \left( \prod_{i=0}^{m-2}
(1+r+(k-1) +(k-1)i) \right)-\left(\prod_{i=0}^{m-2} (r+(k-1)+ (k-1)i) \right) \right) =\\
&&\left( \prod_{i=1}^{m-1} (1+r+(k-1)i) \right) + \\
&&r \sum_{p=0}^{m-2} \left( \prod_{i=0}^{p-1} (r+(k-1) + (k-1)i)
\right)\left( \prod_{i=p+1}^{m-2} (1+r+(k-1)+ (k-1)i) \right) =
\end{eqnarray*}

\begin{eqnarray*}
&&\left( \prod_{i=1}^{m-1} (1+r+(k-1)i) \right)+ \\
&& r \sum_{p=0}^{m-2} \left( \prod_{i=1}^{p} (r+(k-1)i) \right)\left( \prod_{i=p+2}^{m-1} (1+r+ (k-1)i) \right) = \\
&&\left( \prod_{i=1}^{m-1} (1+r+(k-1)i) \right) +r
\sum_{q=1}^{m-1} \left( \prod_{i=0}^{q-1} (r+(k-1)i) \right)
\left( \prod_{i=q+1}^{m-1} (1+r+(k-1)i) \right) = \\
&&\sum_{p=0}^{m-1} \left( \prod_{i=0}^{p-1} (r+(k-1)i) \right)
\left( \prod_{i=p+1}^{m-1} (1+r+(k-1)i) \right).
\end{eqnarray*}
Here we have used our induction hypothesis for $\Omega(k,m-2, r+k-1)$ in
going from line 4 to line 5 in the above series of equalities.
\end{proof}

Next we turn to general formulas for $B^{(k)}_{s,kn+j}$,
$B^{(k)}_{0,s,kn+j}$ and $B^{(k)}_{1,s,kn+j}$ for all
$k \geq 2$, $n \geq 0$, and
$0 \leq j \leq k-1$ . First, we develop
formulas for the coefficients $B^{(k)}_{s,kn+j}$ by iterating
the recursions (\ref{rec:Cj}) and (\ref{rec:Ck-1}) starting with
our formula for  $B^{(k)}_{n,kn+j}$ given in   Theorem \ref{thm:B4}.
That is, we shall prove the following theorem.

\begin{theorem}\label{thm:Cform}
For all $n \geq 0$, $k \geq 2$, and $0 \leq j \leq k-1$,
\begin{equation}\label{eq:Cform}
B^{(k)}_{n-s,kn+j} = ((k-1)n+j)! \sum_{r=0}^s (-1)^{s-r}
\binom{(k-1)n +j +r}{r} \binom{kn +j +1}{s-r} \prod_{i=0}^{n-1}
(r+j +(k-1)i).
\end{equation}
\end{theorem}
\begin{proof}
Our formula holds for $s=0$ by Theorem \ref{thm:B4}.

Now assume that (\ref{eq:Cform}) holds for $s$ for all $n$, $k$, and $j$.
Plugging $n-s$ for $s$ in (\ref{rec:Cj}), we obtain the recursion
\begin{equation}\label{Cform1}
B^{(k)}_{n-s,kn+j+1} = (kn+j+1 -s)B^{(k)}_{n-s,kn+j} +
(s+1)B^{(k)}_{n-s -1,kn+j}
\end{equation}
or, equivalently,
\begin{equation}\label{Cform1.1}
(s+1)B^{(k)}_{n-(s+1),kn+j}= B^{(k)}_{n-s,kn+j+1} - (kn+j+1
-s)B^{(k)}_{n-s,kn+j}
\end{equation}
Thus by our induction hypothesis,
\begin{eqnarray*}
&&(s+1)B^{(k)}_{n-(s+1),kn+j}= \\
&&((k-1)n+j+1)!
\sum_{r=0}^s (-1)^{s-r} \binom{(k-1)n +j +1 +r}{r}
\binom{kn +j +2}{s-r} \prod_{i=0}^{n-1} (r+j +1 +(k-1)i)-\\
&&(kn+j+1 -s) ((k-1)n+j)! \sum_{r=0}^s (-1)^{s-r} \binom{(k-1)n +j
+r}{r} \binom{kn +j +1}{s-r} \prod_{i=0}^{n-1} (r+j +(k-1)i).
\end{eqnarray*}
It follows that
\begin{eqnarray*}
&&\frac{B^{(k)}_{n-(s+1),kn+j}}{((k-1)n+j)!}= \\
&&\frac{(k-1)n+j+1}{s+1}
\sum_{r=0}^s (-1)^{s-r} \binom{(k-1)n +j +1+r}{r}
\binom{kn +j +2}{s-r} \prod_{i=0}^{n-1} (r+j +1 +(k-1)i) \\
&& - \frac{kn+j+1 -s}{s+1}
\sum_{r=0}^s (-1)^{s-r} \binom{(k-1)n +j +r}{r}
\binom{kn +j +1}{s-r} \prod_{i=0}^{n-1} (r+j +(k-1)i).
\end{eqnarray*}
We can then look at three terms. First the $r=s$ term on the first
summand above is
\begin{eqnarray}\label{Cform2}
&&\frac{(k-1)n+j+1}{s+1}\binom{(k-1)n +j +1+s}{s} \prod_{i=0}^{n-1} (s+j +1 +(k-1)i) = \nonumber \\
&&\binom{(k-1)n +j +1+s}{s+1}\prod_{i=0}^{n-1} (s+1+j +(k-1)i).
\end{eqnarray}
Second, the $r=0$ term from the second summand is
\begin{equation}\label{Cform3}
(-1)^{s+1} \frac{kn+j+1 -s}{s+1}\binom{kn +j +1}{s}
\prod_{i=0}^{n-1} (j +(k-1)i)= (-1)^{s+1} \binom{kn +j
+1}{s+1}\prod_{i=0}^{n-1} (j +(k-1)i).
\end{equation}
Finally we can combine the terms involving $\prod_{i=0}^{n-1} (r+j +(k-1)i)$ from the first and second summands to obtain
\begin{eqnarray}\label{Cform4}
&& \sum_{r=1}^s (-1)^{s+1-r} \prod_{i=0}^{n-1} (r+j +(k-1)i) \times
\nonumber \\
&&\left[ \frac{(k-1)n+j+1}{s+1} \binom{(k-1)n +j +1 +r-1}{r-1}
\binom{kn +j +2}{s-(r-1)} \right.\nonumber \\
&& \left. + \frac{kn+j+1 -s}{s+1} \binom{(k-1)n +j +r}{r}
\binom{kn +j +1}{s-r} \right]
\end{eqnarray}
If we let $D_r$ be the factor in square brackets in the $r$-th term
of (\ref{Cform4}), then we see that
\begin{eqnarray}\label{Cform5}
D_r &=& \frac{r}{r}\cdot\frac{(k-1)n+j+1}{s+1} \binom{(k-1)n +j
+r}{r-1}
\binom{kn+j+2}{s+1-r}\nonumber +\\
&& \frac{s+1-r}{s+1-r}\cdot\frac{kn+j+1 -s}{s+1} \binom{(k-1)n +j
+r}{r} \binom{kn +j +1}{s-r} \nonumber \end{eqnarray}
\begin{eqnarray}
&=&\frac{r}{s+1} \binom{(k-1)n +j +r}{r}
\binom{kn+j+2}{s+1-r} + \nonumber \\
&& \frac{s+1-r}{s+1-r}\cdot\frac{kn+j+1 -s}{s+1} \binom{(k-1)n +j
+r}{r} \binom{kn +j +1}{s-r} \nonumber \end{eqnarray}
\begin{eqnarray}&=& \binom{(k-1)n +j +r}{r}\left[ \frac{r}{s+1}
\binom{kn+j+2}{s+1-r} + \frac{s+1-r}{s+1-r}\cdot\frac{kn+j+1
-s}{s+1}\binom{kn +j +1}{s-r} \right] \nonumber \end{eqnarray}
\begin{eqnarray}&=& \binom{(k-1)n +j +r}{r}
\frac{(kn+j+1)\downarrow_{s-r}}{(s+1-r)!} \frac{1}{s+1}
\left[r(kn+j+2) +(s+1-r)(kn+j+1-s)\right] \nonumber \end{eqnarray}
\begin{eqnarray} &=& \binom{(k-1)n +j +r}{r}
\frac{(kn+j+1)\downarrow_{s-r}}{(s+1-r)!}
\frac{1}{s+1} \left[ (s+1) (kn+j+1-(s-r)) \right] \nonumber \\
&=&  \binom{(k-1)n +j +r}{r} \frac{(kn+j+1)\downarrow_{s+1-r}}{(s+1-r)!}
\nonumber \\
&=& \binom{(k-1)n +j +r}{r} \binom{kn +j +1}{s+1-r}.
\end{eqnarray}
Thus plugging in (\ref{Cform5}) into (\ref{Cform4}), we see that
(\ref{Cform4}) is equal to
\begin{equation}\label{Cform6}
\sum_{r=1}^s (-1)^{s+1-r} \binom{(k-1)n +j +r}{r} \binom{kn +j +1}{s+1-r}
\prod_{i=0}^{n-1} (r+j +(k-1)i).
\end{equation}
Finally, combining (\ref{Cform2}), (\ref{Cform3}), and
(\ref{Cform6}), we obtain that
\begin{equation}\label{Cform7}
\frac{B^{(k)}_{n-(s+1),kn+j}}{((k-1)n+j)!}= \sum_{r=0}^{s+1}
(-1)^{s+1-r} \binom{(k-1)n +j +r}{r} \binom{kn +j +1}{s+1-r}
\prod_{i=0}^{n-1} (r+j +(k-1)i)
\end{equation}
as desired.
Thus we have established our formula for all $n \geq 0$, $k \geq 2$,
$j =0, \ldots, k-2$ for $s+1$.

Next we substitute  $n-s$ for $s$ in (\ref{rec:Ck-1}), we obtain the recursion
that for $0 \leq j \leq k-2$,
\begin{equation}\label{Cform8}
B^{(k)}_{n-s,kn+k} = (kn+k -s)B^{(k)}_{n-s,kn+k-1} +
(s+1)B^{(k)}_{n-s -1,kn+k-1}
\end{equation}
or, equivalently,
\begin{equation}\label{Cform9}
(s+1)B^{(k)}_{n-(s+1),kn+k-1}= B^{(k)}_{(n+1) -(s+1),k(n+1)} -
(kn+k -s)B^{(k)}_{n-s,kn+k-1}.
\end{equation}
Given that we have already established our formula for $B^{(k)}_{n
-s,kn+k-1}$ and using our induction hypothesis, we obtain that
\begin{eqnarray*}
&&(s+1)B^{(k)}_{n-(s+1),kn+k-1}= \\
&&((k-1)(n+1))! \times \\
&&\sum_{r=0}^{s+1} (-1)^{s+1-r} \binom{(k-1)(n+1) +r}{r}
\binom{k(n+1) +1}{s+1-r} \prod_{i=0}^{n} (r +(k-1)i) -\\
&&(kn+k -s) ((k-1)n+k-1)! \times \\
&&\sum_{r=0}^s (-1)^{s-r} \binom{(k-1)n +k-1 +r}{r}
\binom{kn +k-1 +1}{s-r} \prod_{i=0}^{n-1} (r+(k-1) +(k-1)i).
\end{eqnarray*}
It follows that
\begin{eqnarray*}
&&\frac{B^{(k)}_{n-(s+1),kn+k-1}}{((k-1)n+k-1)!}= \\
&&\frac{1}{s+1} \sum_{r=0}^{s+1} (-1)^{s+1-r}
\binom{(k-1)n+k-1+r}{r}
\binom{kn +k+1}{s+1-r} r \prod_{i=1}^{n} (r+(k-1)i)- \\
&& \frac{kn+k -s}{s+1} \sum_{r=0}^s (-1)^{s-r} \binom{(k-1)n +k-1
+r}{r} \binom{kn +k}{s-r} \prod_{i=0}^{n-1} (r+(k-1) +(k-1)i).
\end{eqnarray*}
We can then look at three terms. First the $r=s+1$ term on the first
summand of the final expression above is
\begin{eqnarray}\label{Cform10}
&&\frac{1}{s+1}\binom{(k-1)n +k-1 +s+1}{s+1} (s+1)
\prod_{i=0}^{n-1}
(s+1+(k-1) +(k-1)i) = \nonumber \\
&&\binom{(k-1)n +k-1+s+1}{s+1}\prod_{i=0}^{n-1} (s+1+(k-1) +(k-1)i).
\end{eqnarray}
Second, the $r=0$ term from the second summand of the final
expression for $\frac{B^{(k)}_{n-(s+1),kn+k-1}}{((k-1)n+k-1)!}$
above is
\begin{equation}\label{Cform11}
(-1)^{s+1} \frac{kn+k -s}{s+1}\binom{kn +k}{s} \prod_{i=0}^{n-1}
((k-1) +(k-1)i)= (-1)^{s+1}\binom{kn +k}{s+1}\prod_{i=0}^{n-1}
((k-1) +(k-1)i).
\end{equation}
Finally we can combine the terms involving $\prod_{i=0}^{n-1}
(r+(k-1) +(k-1)i)$ from the first and second summands of the final
expression for $\frac{B^{(k)}_{n-(s+1),kn+k-1}}{((k-1)n+k-1)!}$
above to obtain
\begin{eqnarray}\label{Cform12}
&& \sum_{r=1}^s (-1)^{s+1-r} \binom{(k-1)n+k-1+r}{r}
\prod_{i=0}^{n-1} (r+(k-1) +(k-1)i) \times \nonumber \\
&&\left[ \frac{r}{s+1} \binom{kn +k +1}{s+1-r} + \frac{kn+k -s}{s+1}  \binom{kn +k}{s-r} \right]
\end{eqnarray}
If we let $E_r$ be the factor in square brackets in the $r$-th term
of (\ref{Cform12}), then we see that
\begin{eqnarray}\label{Cform13}
E_r &=& \frac{r}{s+1} \binom{kn +k +1}{s+1-r}
+ \frac{s+1-r}{s+1-r}\cdot\frac{kn+k -s}{s+1}  \binom{kn +k}{s-r} \nonumber \\
&=& \frac{(kn+k)\downarrow_{s-r}}{(s+1-r)!} \cdot\frac{1}{s+1}
\left[r(kn+k+1) + (kn+k-s)(s+1-r)\right] \nonumber \end{eqnarray}
\begin{eqnarray}&=& \frac{(kn+k)\downarrow_{s-r}}{(s+1-r)!}\cdot \frac{1}{s+1}
\left[(s+1)(kn+k-(s-r))\right] \nonumber \end{eqnarray}
\begin{eqnarray} &=& \frac{(kn+k)\downarrow_{s+1-r}}{(s+1-r)!} =
\binom{kn+k}{s+1-r}.
\end{eqnarray}
Thus plugging in (\ref{Cform13}) into (\ref{Cform12}),
we see that (\ref{Cform12}) is equal to
\begin{equation}\label{Cform14}
\sum_{r=1}^s (-1)^{s+1-r} \binom{(k-1)n +k-1 +r}{r} \binom{kn
+k}{s+1-r} \prod_{i=0}^{n-1} (r+(k-1) +(k-1)i).
\end{equation}
Finally, combining (\ref{Cform10}), (\ref{Cform11}), and (\ref{Cform14}),
we obtain that
\begin{eqnarray*}
&&\frac{B^{(k)}_{n-(s+1),kn+k-1}}{((k-1)n+k-1)!}=\\
&&\sum_{r=0}^{s+1} (-1)^{s+1-r} \binom{(k-1)n +k-1 +r}{r}
\binom{kn +k}{s+1-r} \prod_{i=0}^{n-1} (r+(k-1) +(k-1)i)
\end{eqnarray*}
as desired.
\end{proof}

As a corollary to Theorem~\ref{thm:Cform} we get combinatorial
proofs for two special cases of the Saalsch\"utz's identity. Recall
that in terms of generalized hypergeometric functions can be written
as
$${}_3F_2\left[\begin{array}{ccc} a & b & c \\ d & e & \end{array}; 1\right]=
\frac{(d-a)_{|c|}(d-b)_{|c|}}{d_{|c|}(d-a-b)_{|c|}}$$ where
$d+e=a+b+c+1$ and $c$ is a negative integer\footnote{For the first
identity in Corollary~\ref{saalschutz}, $a=n+2$, $b=n+1$, $c=s-n$,
$d=1$, and $e=n+s+3$; for the second identity there, $a=n+1$, $b=n$,
$c=s-n$, $d=0$, and $e=n+s+2$}.

\begin{corollary}\label{saalschutz} The following identities hold:
$$\frac{n+1}{s+1}\binom{n}{s}^2=\sum_{r=0}^{n-s}(-1)^{n-s-r}\binom{n+r}{r}\binom{n+r+1}{r}\binom{2n+2}{n-s-r};$$
$$\binom{n-1}{s}\binom{n+1}{s+1}=\sum_{r=0}^{n-s}(-1)^{n-s-r}\binom{n+r}{r}\binom{n+r-1}{r-1}\binom{2n+1}{n-s-r}.$$
\end{corollary}

\begin{proof} The RHS of the first identity is $B^{(2)}_{s,2n+1}/(n!(n+1)!)$ (we use Theorem~\ref{thm:Cform}
for $k=2$, $j=1$ and $s=n-s$). However, $B^{(2)}_{s,2n+1}$ counts
exactly the same objects as
$P_{s,2n+1}=P_{0,s,2n+1}+P_{1,s,2n+1}=\frac{1}{s+1}\binom{n}{s}^2((n+1)!)^2$
does in~\cite{kitrem}. The RHS of the second identity is
$B^{(2)}_{s,2n}/(n!)^2$, and $B^{(2)}_{s,2n}$ is the same as
$P_{s,2n}=P_{0,s,2n}+P_{1,s,2n}=(n!)^2\binom{n-1}{s}\binom{n+1}{s+1}$
in~\cite{kitrem}.
\end{proof}

Next we prove formulas for the coefficients $B^{(k)}_{1,s,kn+j}$ by
iterating the recursions (\ref{rec:B1j}) and (\ref{rec:B1k-1})
starting from our formula for $B^{(k)}_{1,n-1,kn+j}$ given in
Theorem \ref{thm:B4:1}.

\begin{theorem}\label{thm:B1form}
For all $n \geq 0$, $k \geq 2$, and $0 \leq j \leq k-1$,
\begin{equation}\label{eq:Bform}
B^{(k)}_{1,n-1-s,kn+j} = ((k-1)n+j)!
\sum_{r=0}^s (-1)^{s-r} \binom{(k-1)n +j +r}{r}
\binom{kn +j }{s-r} \Omega(k,n,r+j).
\end{equation}
where $\Omega(k,n,r)$ is given by~\rm{(}\ref{omega1}\rm{)}.
\end{theorem}

\begin{proof}
Our formula holds for $s=0$ by Theorem \ref{thm:B4:1}.

Now assume that(\ref{eq:Bform}) holds for $s$ for all $n$, $k$, and $j$.
Plugging $n-s-1$ for $s$ in (\ref{rec:B1j}), we obtain the recursion
that for $0 \leq j \leq k-2$,
\begin{equation}\label{Bform1}
B^{(k)}_{1,n-1-s,kn+j+1} = (kn+j -s)B^{(k)}_{1,n-1-s,kn+j} + (s+1)B^{(k)}_{1,n-1-s-1,kn+j}
\end{equation}
or, equivalently,
\begin{equation}\label{Bform1.1}
(s+1)B^{(k)}_{1,n-1-(s+1),kn+j}= B^{(k)}_{1,n-1-s,kn+j+1} -
(kn+j-s)B^{(k)}_{1,n-1-s,kn+j}.
\end{equation}
Thus by our induction hypothesis,
\begin{eqnarray*}
&&(s+1)B^{(k)}_{1,n-1-(s+1),kn+j}= \\
&&((k-1)n+j+1)!
\sum_{r=0}^s (-1)^{s-r} \binom{(k-1)n +j +1 +r}{r}
\binom{kn +j +1}{s-r} \Omega(k,n,r+j+1) \\
&& - (kn+j -s) ((k-1)n+j)!
\sum_{r=0}^s (-1)^{s-r} \binom{(k-1)n +j +r}{r}
\binom{kn +j}{s-r} \Omega(k,n,r+j).
\end{eqnarray*}
It follows that
\begin{eqnarray*}
&&\frac{B^{(k)}_{1,n-1-(s+1),kn+j}}{((k-1)n+j)!}= \\
&&\frac{(k-1)n+j+1}{s+1}
\sum_{r=0}^s (-1)^{s-r} \binom{(k-1)n +j +1+r}{r}
\binom{kn +j +1}{s-r} \Omega(k,n,r+j+1) \\
&& - \frac{kn+j -s}{s+1}
\sum_{r=0}^s (-1)^{s-r} \binom{(k-1)n +j +r}{r}
\binom{kn +j}{s-r} \Omega(k,n,r+j).
\end{eqnarray*}
We can then look at three terms. First the $r=s$ term on the first
summand above is
\begin{eqnarray}\label{Bform2}
&&\frac{(k-1)n+j+1}{s+1}\binom{(k-1)n +j +1+s}{s}
\Omega(k,n,s+j+1) =
\nonumber \\
&&\binom{(k-1)n +j +1+s}{s+1}\Omega(k,n,s+j+1).
\end{eqnarray}
Second, the $r=0$ term from the second summand is
\begin{equation}\label{Bform3}
(-1)^{s+1} \frac{kn+j -s}{s+1}\binom{kn +j}{s} \Omega(k,n,j)=
(-1)^{s+1}\binom{kn +j}{s+1}\Omega(k,n,j).
\end{equation}
Finally we can combine the terms involving $\Omega(k,n,r+j)$ from the first and second summands to obtain
\begin{eqnarray}\label{Bform4}
&& \sum_{r=1}^s (-1)^{s+1-r} \Omega(k,n,r+j) \times \nonumber \\
&&\left[ \frac{(k-1)n+j+1}{s+1} \binom{(k-1)n +j +1 +r-1}{r-1}
\binom{kn +j +1}{s-(r-1)} \right. \nonumber \\
&& \left. + \frac{kn+j -s}{s+1} \binom{(k-1)n +j +r}{r}
\binom{kn +j}{s-r} \right]
\end{eqnarray}
If we let $F_r$ be the factor in square brackets in the $r$-th term
of (\ref{Bform4}), then we see that
\begin{eqnarray}\label{Bform5}
F_r &=& \frac{r}{r}\cdot\frac{(k-1)n+j+1}{s+1} \binom{(k-1)n +j
+r}{r-1}
\binom{kn+j+1}{s+1-r} \nonumber \\
&& + \frac{s+1-r}{s+1-r}\cdot\frac{kn+j -s}{s+1} \binom{(k-1)n +j
+r}{r}
\binom{kn +j}{s-r} \nonumber \\
&=&\frac{r}{s+1} \binom{(k-1)n +j +r}{r}
\binom{kn+j+1}{s+1-r}+ \nonumber \\
&& \frac{s+1-r}{s+1-r}\cdot\frac{kn+j -s}{s+1} \binom{(k-1)n +j
+r}{r}
\binom{kn +j}{s-r} \nonumber \\
&=& \binom{(k-1)n +j +r}{r}\left[ \frac{r}{s+1}
\binom{kn+j+1}{s+1-r} + \frac{s+1-r}{s+1-r}\cdot\frac{kn+j
-s}{s+1}\binom{kn +j}{s-r} \right]
\nonumber \\
&=& \binom{(k-1)n +j +r}{r} \frac{(kn+j)\downarrow_{s-r}}{(s+1-r)!}
\cdot\frac{1}{s+1} \left[r(kn+j+1) +(s+1-r)(kn+j-s)\right] \nonumber \\
&=& \binom{(k-1)n +j +r}{r}
\frac{(kn+j)\downarrow_{s-r}}{(s+1-r)!}
\cdot\frac{1}{s+1} \left[ (s+1) (kn+j-(s-r)) \right] \nonumber \\
&=&  \binom{(k-1)n +j +r}{r}
\frac{(kn+j)\downarrow_{s+1-r}}{(s+1-r)!}
\nonumber \\
&=& \binom{(k-1)n +j +r}{r} \binom{kn +j}{s+1-r}.
\end{eqnarray}
Thus plugging in (\ref{Bform5}) into (\ref{Bform4}),
we see that (\ref{Bform4}) is equal to
\begin{equation}\label{Bform6}
\sum_{r=1}^s (-1)^{s+1-r} \binom{(k-1)n +j +r}{r} \binom{kn +j}{s+1-r}
\Omega(k,n,r+j).
\end{equation}
Finally, combining (\ref{Bform2}), (\ref{Bform3}), and
(\ref{Bform6}), we obtain that
\begin{equation}\label{Bform7}
\frac{B^{(k)}_{1,n-1-(s+1),kn+j}}{((k-1)n+j)!}=
\sum_{r=0}^{s+1} (-1)^{s+1-r} \binom{(k-1)n +j +r}{r} \binom{kn +j}{s+1-r}
\Omega(k,n,r+j)
\end{equation}
as desired.
Thus we have established our formula for all $n \geq 0$, $k \geq 2$,
$j =0, \ldots, k-2$ for $s+1$.

Next we substitute  $n-s-1$ for $s$ in (\ref{rec:B1k-1}), we obtain the
recursion
\begin{equation}\label{Bform8}
B^{(k)}_{1,n-1-s,kn+k} = (kn+k-1 -s)B^{(k)}_{1,n-1-s,kn+k-1} + (s+2)B^{(k)}_{1,n-1-s -1,kn+k-1}
+B^{(k)}_{0,n-1-s,kn+k-1}
\end{equation}
or, equivalently,
\begin{eqnarray}\label{Bform9}
&&(s+2)B^{(k)}_{1,n-1-(s+1),kn+k-1} = \\
&& B^{(k)}_{1,(n+1)-1 -(s+1),k(n+1)} -
(kn+k-1 -s)B^{(k)}_{1,n-1-s,kn+k-1} - B^{(k)}_{0,n-1-s,kn+k-1} = \nonumber \\
&& B^{(k)}_{1,(n+1)-1 -(s+1),k(n+1)} - (kn+k
-1-s)B^{(k)}_{1,n-1-s,kn+k-1}-
\nonumber \\
&&
(B^{(k)}_{n-(s+1),kn+k-1} -B^{(k)}_{1,n-1-s,kn+k-1}) = \nonumber \\
&&B^{(k)}_{1,(n+1)-1 -(s+1),k(n+1)} -
(kn+k-s-2)B^{(k)}_{1,n-1-s,kn+k-1} - B^{(k)}_{n-(s+1),kn+k-1}.
\nonumber
\end{eqnarray}

Given that we have already established our formula for $B^{(k)}_{n
-(s+1),kn+k-1}$ and $B^{(k)}_{1,n-1 -s,kn+k-1}$, and using our
induction hypothesis, we obtain that
\begin{eqnarray*}
&&(s+2)B^{(k)}_{1,n-1-(s+1),kn+k-1}= \\
&&((k-1)(n+1))!
\sum_{r=0}^{s+1} (-1)^{s+1-r} \binom{(k-1)(n+1) +r}{r}
\binom{k(n+1)}{s+1-r} \Omega(k,n+1,r) \\
&& - (kn+k -s-2) ((k-1)n+k-1)! \times \\
&&\sum_{r=0}^s (-1)^{s-r} \binom{(k-1)n +k-1 +r}{r}
\binom{kn +k-1 }{s-r} \Omega(k,n,r+k-1) \\
&& - ((k-1)n+k-1)!\times \\
&&\sum_{r=0}^{s+1} (-1)^{s+1-r} \binom{(k-1)(n+1) +r}{r} \binom{kn
+k}{s+1-r} \prod_{i=0}^{n-1} (r +(k-1) + (k-1)i).
\end{eqnarray*}
It follows that
\begin{eqnarray*}
&&\frac{B^{(k)}_{1,n-1-(s+1),kn+k-1}}{((k-1)n+k-1)!}= \\
&&\frac{1}{s+2} \sum_{r=0}^{s+1} (-1)^{s+1-r}
\binom{(k-1)n+k-1+r}{r}
\binom{kn +k}{s+1-r} \times \\
&& [ \Omega(k,n+1,r) -
\prod_{i=0}^{n-1} (r +(k-1) + (k-1)i)] \\
&& - \frac{kn+k -s-2}{s+2}
\sum_{r=0}^s (-1)^{s-r} \binom{(k-1)n +k-1 +r}{r}
\binom{kn +k-1}{s-r} \Omega(k,n,r+k-1).
\end{eqnarray*}
Next we shall rewrite the term $[ \Omega(k,n+1,r) -
\prod_{i=0}^{n-1} (r +(k-1) + (k-1)i)]$. Note that
\begin{eqnarray}\label{omega}
&&\Omega(k,n+1,r) -
\prod_{i=0}^{n-1} (r +(k-1) + (k-1)i) = \\
&&\sum_{p=0}^n \left( \prod_{i=0}^{p-1} (r+(k-1)i) \right)
\left(\prod_{i=p+1}^{n} (1+r+(k-1)i)  \right) \nonumber \\
&&- \prod_{i=0}^{n-1} (r +(k-1) + (k-1)i) = \nonumber  \\
&&\sum_{p=1}^n r \left( \prod_{i=1}^{p-1} (r+(k-1)i) \right)
\left(\prod_{i=p+1}^{n} (1+r+(k-1)i)  \right) \nonumber \\
&& + \prod_{i=1}^n (1+r+(k-1)i) -
\prod_{i=0}^{n-1} (r +(k-1) + (k-1)i) = \nonumber \\
&&r \sum_{p=0}^{n-1}  \left( \prod_{i=0}^{p-1} (r+(k-1)+(k-1)i)
\right)
\left(\prod_{i=p+1}^{n-1} (1+r+(k-1)+(k-1)i)  \right)\nonumber \\
&& + \prod_{i=0}^{n-1} (1+r+(k-1)+(k-1)i) -
\prod_{i=0}^{n-1} (r +(k-1) + (k-1)i) = \nonumber \\
&& r\Omega(k,n,r+k-1) + \prod_{i=0}^{n-1} (1+r+(k-1)+(k-1)i) -
\prod_{i=0}^{n-1} (r +(k-1) + (k-1)i) = \nonumber \\
&& r\Omega(k,n,r+k-1) + \Omega(k,n,r+k-1) \ =\
(1+r)\Omega(k,n,r+k-1). \nonumber
\end{eqnarray}
Here on the first to last line, we have used Lemma
\ref{lem:omega}. Thus plugging~(\ref{omega}) into our expression
for $\frac{B^{(k)}_{1,n-1-(s+1),kn+k-1}}{((k-1)n+k-1)!}$, we
obtain that
\begin{eqnarray*}
&&\frac{B^{(k)}_{1,n-1-(s+1),kn+k-1}}{((k-1)n+k-1)!} = \\
&&\frac{1}{s+2} \sum_{r=0}^{s+1} (-1)^{s+1-r}
\binom{(k-1)n+k-1+r}{r}
\binom{kn +k}{s+1-r}  (1+r)\Omega(k,n,r+k-1)\\
&& - \frac{kn+k -s-2}{s+2}
\sum_{r=0}^s (-1)^{s-r} \binom{(k-1)n +k-1 +r}{r}
\binom{kn +k-1}{s-r} \Omega(k,n,r+k-1).
\end{eqnarray*}
We can then look at two terms. First the $r=s+1$ term on the first
summand of the final expression above is
\begin{eqnarray}\label{Bform10}
&&\frac{1}{s+2}\binom{(k-1)n +k-1 +s+1}{s+1} (s+2)
\Omega(k,n,s+1 +k-1) = \nonumber \\
&&\binom{(k-1)n +k+s}{s+1}\Omega(k,n,s+k).
\end{eqnarray}
Second,  we can combine the terms involving $\Omega(k,n,r +k-1)$ from the first and second summands of the final expression for
$\frac{B^{(k)}_{1,n-1-(s+1),kn+k-1}}{((k-1)n+k-1)!}$ above to obtain
\begin{eqnarray}\label{Bform12}
&& \sum_{r=0}^s (-1)^{s+1-r} \binom{(k-1)n+k-1+r}{r}
\Omega(k,n,r +k-1) \times \nonumber \\
&&\left[ \frac{1+r}{s+2} \binom{kn +k}{s+1-r}
+ \frac{kn+k -s-2}{s+2}  \binom{kn +k-1}{s-r} \right]
\end{eqnarray}
If we let $G_r$ be the factor in square brackets in the $r$-th term
of (\ref{Bform12}), then we see that
\begin{eqnarray}\label{Bform13}
G_r &=& \frac{1+r}{s+2} \binom{kn +k}{s+1-r}
+ \frac{s+1-r}{s+1-r}\cdot\frac{kn+k -s-2}{s+2}  \binom{kn +k-1}{s-r} \nonumber \\
&=& \frac{(kn+k-1)\downarrow_{s-r}}{(s+1-r)!} \frac{1}{s+2}
\left[(1+r)(kn+k) + (kn+k-(s+2))(s+1-r)\right] \nonumber \\
&=& \frac{(kn+k-1)\downarrow_{s-r}}{(s+1-r)!} \frac{1}{s+2}
\left[(s+2)(kn+k-1-(s-r))\right] \nonumber \\
&=& \frac{(kn+k-1)\downarrow_{s+1-r}}{(s+1-r)!} \nonumber \\
&=& \binom{kn+k-1}{s+1-r}.
\end{eqnarray}
Thus plugging in (\ref{Bform13}) into (\ref{Bform12}),
we see that (\ref{Bform12}) is equal to
\begin{equation}\label{Bform14}
\sum_{r=0}^s (-1)^{s+1-r} \binom{(k-1)n +k-1 +r}{r} \binom{kn +k-1}{s+1-r}
\Omega(k,n,r+k-1).
\end{equation}
Finally, combining (\ref{Bform10}) and (\ref{Bform14}), we obtain
that
\begin{equation}\label{Bform15}
\frac{B^{(k)}_{1,n-1-(s+1),kn+k-1}}{((k-1)n+k-1)!}=
\sum_{r=0}^{s+1} (-1)^{s+1-r} \binom{(k-1)n +k-1 +r}{r} \binom{kn +k-1}{s+1-r} \Omega(k,n,r+k-1)
\end{equation}
as desired.
\end{proof}

Note that it now follows that for all $n \geq 0$, $k \geq 2$, $0
\leq j \leq k-1$, and $0 \leq s \leq n-1$,\\

$B^{(k)}_{0,n-1-s,kn+j} = B^{(k)}_{n-(s+1),kn+j}-
B^{(k)}_{1,n-1-s,kn+j}=$
\begin{eqnarray*} &&
((k-1)n+j)! \sum_{r=0}^{s+1} (-1)^{s+1-r} \binom{(k-1)n +j  +r}{r}
\binom{kn +j+1}{s+1-r} \prod_{i=0}^{n-1} (r+j+(k-1)i)\\
&&  - {((k-1)n+j)!} \sum_{r=0}^{s} (-1)^{s-r} \binom{(k-1)n +j
+r}{r} \binom{kn +j}{s-r} \Omega(k,n,r+j).
\end{eqnarray*}
Thus
\begin{eqnarray*}
&&\frac{B^{(k)}_{0,n-(s+1),kn+j}}{((k-1)n+j)!}= \\
&&\binom{(k-1)n +j  +s+1}{s+1} \prod_{i=0}^{n-1} (s+1+j+(k-1)i)\\
&& +
\sum_{r=0}^{s} (-1)^{s+1-r} \binom{(k-1)n +j  +r}{r} \binom{kn +j+1}{s+1-r} \prod_{i=0}^{n-1} (r+j+(k-1)i)\\
&&  -
\sum_{r=0}^{s} (-1)^{s-r} \binom{(k-1)n +j +r}{r} \binom{kn +j}{s-r}
\times \\
&&\left[\prod_{i=0}^{n-1} (1+r +j +(k-1)i) - \prod_{i=0}^{n-1} (r +j
+(k-1)i)\right] = \\
&&\binom{(k-1)n +j  +s+1}{s+1} \prod_{i=0}^{n-1} (s+1+j+(k-1)i) \\
&&+ \sum_{r=0}^{s} (-1)^{s+1-r} \binom{(k-1)n +j  +r}{r}
\left[\binom{kn +j+1}{s+1-r} -\binom{kn +j}{s-r}\right]
\prod_{i=0}^{n-1} (r+j+(k-1)i) \\
&&  - \sum_{r=0}^{s} (-1)^{s-r} \binom{(k-1)n +j +r}{r} \binom{kn
+j}{s-r} \prod_{i=0}^{n-1} (1+r +j +(k-1)i) =\end{eqnarray*}

\begin{eqnarray*}
&&\binom{(k-1)n +j  +s+1}{s+1} \prod_{i=0}^{n-1} (s+1+j+(k-1)i)\\
&& +
\sum_{r=0}^{s} (-1)^{s+1-r} \binom{(k-1)n +j  +r}{r}
\binom{kn+j}{s+1-r}\prod_{i=0}^{n-1} (r+j+(k-1)i) \\
&&  -
\sum_{r=0}^{s} (-1)^{s-r} \binom{(k-1)n +k-1 +r}{r} \binom{kn +j}{s-r}
\prod_{i=0}^{n-1} (1+r +j +(k-1)i)=\\
&&
\sum_{r=0}^{s+1} (-1)^{s+1-r} \binom{(k-1)n +j  +r}{r} \binom{kn +j}{s+1-r} \prod_{i=0}^{n-1} (r+j+(k-1)i) \\
&&  -
\sum_{r=0}^{s} (-1)^{s-r} \binom{(k-1)n +k-1 +r}{r} \binom{kn +j}{s-r}
\prod_{i=0}^{n-1} (1+r +j +(k-1)i).
\end{eqnarray*}
Thus we have proved the following.

\begin{theorem}\label{thm:B0}
For all $n \geq 0$, $k \geq 2$, and $0 \leq j \leq k-1$,
\begin{enumerate}
\item $B^{(k)}_{0,n,kn+j} = ((k-1)n+j)! \prod_{i=0}^{n-1}
(j+(k-1)i)$ and

\item $B^{(k)}_{0,n-1-s,kn+j} =$\\
$ \sum_{r=0}^{s+1} (-1)^{s+1-r} \binom{(k-1)n +j  +r}{r} \binom{kn +j}{s+1-r} \prod_{i=0}^{n-1} (r+j+(k-1)i)$\\
$  -
\sum_{r=0}^{s} (-1)^{s-r} \binom{(k-1)n +k-1 +r}{r} \binom{kn +j}{s-r}
\prod_{i=0}^{n-1} (1+r +j +(k-1)i)$ for $0 \leq s \leq n-1$.
\end{enumerate}
\end{theorem}

Finally, we shall end this section by showing that we can get
another set of formulas for the coefficients $B^{(k)}_{s,kn+j}$ for all
$k \geq 2$, $n \geq 0$ by iterating the recursions (\ref{rec:Cj}) and
(\ref{rec:Ck-1}) starting with our formula for $B^{(k)}_{s,kn+j}$ given
in Theorem \ref{thm:B3}.

\begin{theorem}\label{thm:BF}
For all $k \geq 2$, $n \geq 0$, $0 \leq j \leq k-1$, and $0 \leq s
\leq n$,
\begin{eqnarray}\label{eq:BFs,kn+j}
&&B^{(k)}_{s,kn+j} = \\
&&((k-1)n+j)!\left[\sum_{r=0}^{s} (-1)^{s-r} \binom{(k-1)n+j +r}{r}
\binom{kn+j+1}{s-r} \prod_{i=1}^n (1+r+(k-1)i)\right]. \nonumber
\end{eqnarray}
\end{theorem}
\begin{proof}

Again we proceed by induction on $s$. By Theorem~\ref{thm:B3}, we
have proved our formula for $B^{(k)}_{s,kn+j}$
in the case where $s = 0$ for
all $n \geq 0$ and $0 \leq j \leq k-1$.

Now assume that $s > 0$ and that the theorem holds for all $s' < s$
by induction. Note that for $n = 0$ and $j = 0, \ldots, k-1$, our
formula asserts that
\begin{eqnarray*}
B^{(k)}_{s,j} &=& j! \sum_{r=0}^{s} (-1)^{s-r}
\binom{j +r}{r} \binom{j+1}{s-r}\\
&=& j! \left( \sum_{n \geq 0} \binom{j+n}{n} x^n \right) \left(\sum_{ m\geq 0}
(-1)^m \binom{j+1}{m}x^m\right)|_{x^s}\\
&=& j! \frac{1}{(1-x)^{j+1}} (1-x)^{j+1}|_{x^s} = 0\\
\end{eqnarray*}
so that our formula holds for $n =0$ and for $j =0, \ldots, k-1$.

Next, by induction, assume that our formula holds for $s$, for $n' <
n$ and $j= 0, \ldots, k-1$. Recall by (\ref{rec:Ck-1}), we have
\begin{equation}\label{BF1:s}
B^{(k)}_{s,kn} =  (1+s +(k-1)n)B^{(k)}_{s,kn-1} +(n-s)B^{(k)}_{s-1,kn-1}.
\end{equation}
Thus
by induction, we get that
\begin{eqnarray}\label{BF2:s}
&&\frac{B^{(k)}_{s,kn}}{((k-1)(n-1) + k-1)!}  =
\frac{B^{(k)}_{s,kn}}{((k-1)n)!} = \nonumber \\
&&(1+s +(k-1)n) \times \nonumber \\
&&\sum_{r=0}^{s} (-1)^{s-r} \binom{(k-1)(n-1) + k-1 +r}{r} \binom{k(n-1)+k-1 +1}{s-r} \prod_{i=1}^{n-1} (1+r+(k-1)i) \nonumber \\
&& + (n-s) \times \nonumber \\
&&\sum_{r=0}^{s-1} (-1)^{s-1-r} \binom{(k-1)(n-1)+k-1+r}{r}
\binom{k(n-1)+k-1 +1}{s-1-r} \prod_{i=1}^{n-1} (1+r+(k-1)i).
\nonumber
\end{eqnarray}
Thus
\begin{eqnarray}\label{BF3:s}
&&\frac{B^{(k)}_{s,kn}}{((k-1)n)!} = \\
&&(1+s +(k-1)n)\sum_{r=0}^{s} (-1)^{s-r} \binom{(k-1)n+r}{r} \binom{kn}{s-r} \prod_{i=1}^{n-1} (1+r+(k-1)i) \nonumber \\
&& + (n-s) \sum_{r=0}^{s-1} (-1)^{s-1-r} \binom{(k-1)n +r}{r}
\binom{kn}{s-1-r} \prod_{i=1}^{n-1} (1+r+(k-1)i). \nonumber
\end{eqnarray}
Now we can divide LHS of (\ref{BF3:s}) into two parts. First $r=s$ term
in the first summand is
\begin{equation}\label{BF4:s}
(1+s +(k-1)n) \binom{(k-1)n+s}{s} \prod_{i=1}^{n-1} (1+s+(k-1)i) =
\binom{(k-1)n+s}{s} \prod_{i=1}^{n} (1+s+(k-1)i).
\end{equation}
Then we can combine the remaining terms on the LHS of (\ref{BF3:s}) to
obtain
\begin{eqnarray}\label{BF5:s}
&&\sum_{r=0}^{s-1} (-1)^{s-r} \binom{(k-1)n+r}{r}
\prod_{i=1}^n (1+r+(k-1)i) \times \nonumber \\
&& \ \ \frac{1}{1+r+(k-1)n}\left[ ((k-1)n+s+1)\binom{kn}{s-r} -
(n-s) \binom{kn}{s-1-r}\right].
\end{eqnarray}
The term in square brackets in (\ref{BF5:s}) is equal to
\begin{eqnarray}\label{BF6:s}
&&\left[ ((k-1)n+s+1)\binom{kn+1}{s-r}\frac{kn+1-s-r}{kn+1}
\right. \nonumber \\
&&\left.  -(n-s) \binom{kn+1}{s-r} \frac{s-r}{kn+1}
\right] = \nonumber \\
&&\frac{1}{kn+1} \binom{kn+1}{s-r}
[((k-1)n+s+1)(kn+1 -(s-r)) - (n-s)(s-r)] = \nonumber \\
&&\frac{1}{kn+1} \binom{kn+1}{s-r} [(kn+1)(1+r+(k-1)n)] = \nonumber \\
&&(1+r+(k-1)n)\binom{kn+1}{s-r}.
\end{eqnarray}
Thus plugging in (\ref{BF6:s}) into (\ref{BF5:s}), we see that
(\ref{BF5:s}) is equal to
\begin{equation}\label{BF7:s}
\sum_{r=0}^{s-1} (-1)^{s-r} \binom{(k-1)n+r}{r} \binom{kn+1}{s-r}
\prod_{i=1}^n (1+r+(k-1)i).
\end{equation}
Thus we can combine (\ref{BF4:s}) into (\ref{BF7:s}) to obtain that
\begin{equation}\label{BF8:s}
\frac{B^{(k)}_{s,kn}}{((k-1)n)!} = \sum_{r=0}^{s} (-1)^{s-r} \binom{(k-1)n+r}{r} \binom{kn+1}{s-r}\prod_{i=1}^n (1+r+(k-1)i).
\end{equation}
as desired.

By (\ref{rec:Cj}), we have that for
\begin{equation}\label{BF9:s}
B^{(k)}_{s,kn+j+1} =
((k-1)n+j+s+1)B^{(k)}_{s,kn+j} +(n-s+1)B^{(k)}_{s-1,kn+j-1}
\end{equation}
for $0 \leq j \leq k-2$.  We can assume by induction
that our formula holds for $B^{(k)}_{s,kn+j}$ so that
\begin{eqnarray}\label{BF10:s}
&& \frac{B^{(k)}_{s,kn+j+1}}{((k-1)n+j)!} =  \nonumber \\
&& ((k-1)n+j+s+1)\sum_{r=0}^s (-1)^{s-r}\binom{(k-1)n+j+r}{r}
\binom{kn+j+1}{s-r} \prod_{i=1}^n (1+r + (k-1)i) \nonumber \\
&& +(n-s+1)\sum_{r=0}^{s-1} (-1)^{s-1-r}\binom{(k-1)n+j+r}{r}
\binom{kn+j+1}{s-1-r} \prod_{i=1}^n (1+r + (k-1)i) \nonumber.
\end{eqnarray}

Again, we can divide the LHS of (\ref{BF10:s}) into two terms.
First term coming from the $r=s$ of the first summand is
\begin{eqnarray}\label{BF11:s}
&&((k-1)n+j+s+1) \binom{(k-1)n+j+s}{s} \prod_{i=1}^n (s+(k-1)i) = \\
&& ((k-1)n+j+1) \binom{(k-1)n+j+1+s}{s}\prod_{i=1}^n (s+(k-1)i)
\nonumber.
\end{eqnarray}

Next the remaining terms on the LHS of (\ref{BF10:s}) can be
combined into
\begin{eqnarray}\label{BF12:s}
&&\sum_{r=0}^{s-1} (-1)^{s-r} \prod_{i=1}^n (1+r+(k-1)i) \times \nonumber \\
&&\left[ ((k-1)n+j+s+1)\binom{(k-1)n+j+r}{r}
\binom{kn+j+1}{s-r} - \right.\nonumber \\
&&\left. (n-s+1)\binom{(k-1)n+j+r}{r}
\binom{kn+j+1}{s-1-r}\right].
\end{eqnarray}
Using that
\begin{equation}\label{BF13:s}
\binom{(k-1)n+j+r}{r} =
\binom{(k-1)n+j+r+1}{r} \frac{(k-1)n+j+1}{(k-1)n+j+1+r},
\end{equation}
we can rewrite (\ref{BF12:s}) as
\begin{eqnarray}\label{BF14:s}
&&((k-1)n+j+1)
\sum_{r=0}^{s-1} (-1)^{s-r} \binom{(k-1)n+j+r+1}{r}
\prod_{i=1}^n (1+r+(k-1)i) \times \nonumber \\
&&\frac{1}{(k-1)n+j+1+r}\left[ ((k-1)n+j+s+1)\binom{kn+j+1}{s-r} - \right.\nonumber \\
&&\left. (n-s+1)\binom{kn+j+1}{s-1-r}\right].
\end{eqnarray}
Next observe that the term in the square brackets in (\ref{BF14:s}) is
\begin{eqnarray}\label{BF15:s}
&&((k-1)n+j+s+1) \binom{kn+j+2}{s-r}\frac{kn+j+2-(s-r)}{kn+j+2} \nonumber \\
&& - (n-s+1) \binom{kn+j+2}{s-r} \frac{s-r}{kn+j+2} = \nonumber \\
&& \frac{1}{kn+j+2} \binom{kn+j+2}{s-r} \times \nonumber \\
&&\left[ ((k-1)n+j+s) (kn+j+2-(s-r)) - (n-s+1)(s-r)\right]=
\nonumber \\
&& \frac{1}{kn+j+2} \binom{kn+j+2}{s-r} [(kn+j+2)((k-1)n+j+1+r)]=
\nonumber \\
&& ((k-1)n+j+1+r) \binom{kn+j+2}{s-r}.
\end{eqnarray}

Thus substituting in the result of (\ref{BF15:s}) into
(\ref{BF14:s}), we see that (\ref{BF14:s}) is equal to
\begin{equation}\label{BF16:s}
((k-1)n+j+1) \sum_{r=0}^{s-1} (-1)^{s-r} \binom{(k-1)n+j+1+r}{r}
\binom{kn+j+1+1}{s-r}\prod_{i=1}^n (1+r+(k-1)i).
\end{equation}
Combining (\ref{BF11:s}) and (\ref{BF16:s}), we see that
\begin{eqnarray}\label{BF17:s}
&&\frac{B^{(k)}_{s,kn+j+1}}{((k-1)n+j)!} = ((k-1)n+j+1)\times\\
&&\sum_{r=0}^{s} (-1)^{s-r} \binom{(k-1)n+j+1+r}{r}
\binom{kn+j+1+1}{s-r}\prod_{i=1}^n (1+r+(k-1)i) \nonumber
\end{eqnarray}
or, equivalently, that
 \begin{eqnarray}\label{BF18:s}
&&B^{(k)}_{s,kn+j+1} = \nonumber \\
&&((k-1)n+j+1)! \sum_{r=0}^{s} (-1)^{s-r} \binom{(k-1)n+j+1+r}{r}
\binom{kn+j+2}{s-r}\prod_{i=1}^n (1+r+(k-1)i) \nonumber
\end{eqnarray}
as desired.
This completes our induction and hence we have established
our formulas for $B^{(k)}_{s,kn+j}$ for all $k \geq 2$, $n \geq 0$,
$0 \leq j \leq k-1$  and $0 \leq s \leq n$.
\end{proof}

Once again, there are a number of remarkable identities that now
follow. For example, it
follows from Theorems \ref{thm:Cform} and \ref{thm:BF} that

\begin{theorem}\label{thm:Bid1}
For all $n \geq 0$, $k \geq 2$, and $0 \leq j \leq k-1$,
\begin{eqnarray}\label{eq:Bid}
&& \sum_{r=0}^s (-1)^{s-r}
\binom{(k-1)n +j +r}{r} \binom{kn +j +1}{s-r} \prod_{i=1}^{n}
(1+r+(k-1)i) = \nonumber \\
&&\sum_{r=0}^{n-s} (-1)^{n-s-r}
\binom{(k-1)n +j +r}{r} \binom{kn +j +1}{n-s-r} \prod_{i=0}^{n-1}
(r+j +(k-1)i).
\end{eqnarray}
\end{theorem}

\section{Bijective proofs related to the context}\label{section4}

In this section we generalize bijective proofs from~\cite{kitrem}
to show some relations for the coefficients of $A^{(k)}_{n}(x)$
and $B^{(k)}_{n}(x,z)$ for certain $n$. An important observation
here is that the compliment $\pi^c$ for a $(kn+k-1)$-permutation
$\pi$ has the following property: $\pi_i$ is divisible by $k$ if
and only if $\pi^c_i$ is divisible by $k$ for $i=1,2,\ldots,
kn+k-1$.

First it is easy to see that the polynomial $A^{(k)}_{kn+k-1}(x)$
and $B^{(k)}_{kn+k-1}(x,z)$ have simple symmetry properties.
\begin{theorem} For all $k \geq 2$ and $n \geq 0$
\begin{eqnarray*}
A^{(k)}_{s, kn+k-1} &=& A^{(k)}_{n-s, kn+k-1} \ \mbox{for} \
0 \leq s \leq n,\\
B^{(k)}_{0,s, kn+k-1} &=& B^{(k)}_{0,n-s, kn+k-1} \ \mbox{for} \
0 \leq s \leq n, \ \mbox{and} \\
B^{(k)}_{1,s, kn+k-1} &=& B^{(k)}_{1,n-1-s, kn+k-1} \
\mbox{for} \ 0 \leq s \leq n-1.
\end{eqnarray*}
\end{theorem}
\begin{proof}
The map which sends $\sg$ to $\sg^c$ gives a bijective proof of all
three results.
\end{proof}

\begin{lemma}\label{A=B}
For all $n \geq 0$, $k \geq 2$, and $0 \leq s \leq n$,
\begin{equation*}\label{A=Bk-1}
A^{(k)}_{s,kn+k-1} = B^{(k)}_{s,kn+k-1}.
\end{equation*}
\end{lemma}
\begin{proof}
Given $\sg = \sg_1 \cdots \sg_{kn+k-1} \in S_{kn+k-1}$, let $\sg^*
= (kn+k-\sg_{kn+k-1}) \cdots (kn+k-\sg_{1})$ be the permutation
that results by taking the complement of $\sg$  and then the
reversing the result.  It is then easy to see that
$\overleftarrow{des}_{kN}(\sg) = \overrightarrow{des}_{kN}(\sg^*)$
which shows that $A^{(k)}_{s,kn+k-1} = B^{(k)}_{s,kn+k-1}$.
\end{proof}

\begin{theorem}\label{bij01} For all $k \geq 3$,
$n \geq 0$, and $1 \leq j \leq \lfloor n/k \rfloor$,
\begin{equation*}
A^{(k)}_{j,kn+k-2}(x) =  A^{(k)}_{n-j,kn+k-2}(x).
\end{equation*}
Thus $A^{(k)}_{kn+k-2}(x)$ is symmetric for $n \geq 0$ and $k \geq
3$.
\end{theorem}
\begin{proof}

Given a permutation $\sg = \sg_1 \cdots \sg_{kn+k-2} \in
\S_{kn+k-2}$ with $\overleftarrow{des}_E(\sg) =j$, apply the
complement to the permutation $\sg'=\sg (kn+k-1)$, that is, $\sg'$
is obtained from $\sg$ by adding a dummy letter $(kn+k-1)$ at the
end. In the obtained permutation $\sg'^c$, make a cyclic shift to
the left to make the letter $(kn+k-1)$ be the first one. Remove
$(kn+k-1)$ to get a $(kn+k-2)$-permutation $\sg^*$ with
$\overleftarrow{des}_{kN}(\sg^*) =n-j$. To reverse this procedure
adjoin $(kn+k-1)$ from the left to a given permutation $\sg^*$
with $\overleftarrow{des}_{kN}(\sg^*) =n-j$. Then make a cyclic
shift to the right to make $1$ be the rightmost letter. Use the
complement and remove $(kn+k-1)$ from the obtained permutation to
get a permutation $\sg$ with $\overleftarrow{des}_E(\sg) =j$.

The map described above and its reverse are clearly injective. We
only need to justify that given $j$ occurrences of the descents in
$\sg$, we get $(n-j)$ occurrences in $\sg^*$ (the reverse to this
statement will follow using the same arguments). Notice that
adding $(kn+k-1)$ at the end does not increase the number of
descents. Since $\sg'$ ends with a number which is not equivalent
to $k \mod 0$, $\sg'^c$ has $n-j$ descents. If $\sg'=A1B(kn+k-1)$,
where $A$ and $B$ are some factors, then $\sg'^c=A^c(kn+k-1)B^c1$
and $\sg^*$ is $(kn+k-1)B^c1A^c$ without $(kn+k-1)$. The last
thing to observe is that moving $A^c$ to the end of $\sg'^c$ does
not create a new descent since it cannot start with $1$, also we
do not lose any descents since none of them can end with
$(kn+k-1)$. So, $\sg^*$ has $n-k$ descents.

\end{proof}

\begin{theorem} For all $k \geq 3$,
$n \geq 0$, and $1 \leq j \leq \lfloor n/k \rfloor$,
\begin{equation*}
A^{(k)}_{j,kn+k-2} =  B^{(k)}_{0,j,kn+k-2} +
B^{(k)}_{1,j-1,kn+k-2}.
\end{equation*}
\end{theorem}

\begin{proof}
We now keep track of whether or not this is a number divisible by
$k$ to the left of 1 in $\sg= \sg_1\sg_2 \cdots \sg_{kn+k-2} \in
\S_{kn+k-2}$ with $\overleftarrow{des}_{kN}(\sg) =j$. We do not
provide all the justifications in our explanation of the map,
since they are similar to that in the proof of
Theorem~\ref{bij01}.

Suppose $\sg'=\sg (kn+k-1)=Ax1B(kn+k-1)$, where $A$ and $B$ are
some factors and $x$ is a number. Apply the complement {\em and
reverse} to $\sg'$ to get $(\sg')^{cr}=1B^{cr}(kn+k-1)x^cA^{cr}$.
Make a cyclic shift to the left in $(\sg')^{cr}$ to make the
number $(kn+k-1)$ be the first one and to get
$\sg^*=(kn+k-1)x^cA^{cr}1B^{cr}$. One can check that
$\overleftarrow{des}_{kN}(\sg) = \overrightarrow{des}_{kN}(\sg^*)
= j$. Also, $x$ is divisible by $k$ if and only if $x^c$ is
divisible by $k$. Now, if we remove $(kn+k-1)$ from $\sg^*$ and
$x^c$ is 0 mod $k$, we loose one descent obtaining a permutation
counted by $B^{(k)}_{1,j-1,kn+k-2}$; if we remove $(kn+k-1)$ from
$\sg^*$ and $x^c$ is not o mod $k$, the number of descents in the
obtained permutation is the same, $j$, and thus we get a
permutation counted by $B^{(k)}_{0,j,kn+k-2}$. Note that if $Ax$
is the empty word, that is, $\sg$ starts with 1, then this case is
treated as the case ``$x$ is not 0 mod $k$" since $\sg^*$ will
start with $(kn+k-1)1$. Thus one may think of $kn+k-1$ as the
(cyclic) predecessor of 1 in this case, that is, $x=kn+k-1$.

The reverse to the map described is easy to see.
\end{proof}

\section{Open questions}\label{section5}
There are a number of questions that arise from this work. For
example, our proofs of the formulas for $A^{(k)}_{s,kn+j}$,
$B^{(k)}_{s,kn+j}$, $B^{(k)}_{0,s,kn+j}$, and $B^{(k)}_{1,s,kn+j}$
all arise by iterating simple recursion. It would be interesting to
find more direct combinatorial proofs of these facts. For example,
we know that $A^{(2)}_{s,2n} = (n!)^2 \binom{n}{s}^2$. We have
direct combinatorial proofs in the case of $s=0$ and $s=n$. It is
not difficult to extend those proofs to also prove the case $s=1$,
but we do not have a simple combinatorial proof for general $s$.
More generally, we have the following problem.

{\bf Problem 1.} Give an inclusion-exclusion combinatorial argument
to prove the formulas for $A^{(k)}_{s,kn+j}$ and $B^{(k)}_{s,kn+j}$.
In particular, to proceed with $A^{(k)}_{s,kn+j}$, a point to start
might be the case $s=0$ (avoidance) and understanding by
inclusion-exclusion the fact that

$$\begin{array}{l} A_{0,2n}^{(2)}=(n!)^2\sum_{r=0}^{n}(-1)^{n-r}{2n+1
\choose n-r}{n+r \choose n}^2;\\
A_{0,2n+1}^{(2)}=B_{0,2n+1}^{(2)}=(n+1)!n!\sum_{r=0}^{n}(-1)^{n-r}{2n+2
\choose n-r}{n+r \choose n}{n+r+1 \choose r}; \\
B_{0,2n}^{(2)}=(n!)^2\sum_{r=0}^{n}(-1)^{n-r}{2n+1 \choose n-r}{n+r
\choose n}{n+r-1 \choose n}. \end{array}$$

{\bf Problem 2} A more general problem is to study pattern matching
in permutations where one takes into account the equivalences
classes of the element mod $k$ for some $k \geq 2$. That is, given
any sequence $\sg = \sg_1 \cdots \sg_n$ of distinct integers, we let
$red(\sg)$ be the permutation that results by replacing the $i$-th
largest integer that appears in the sequence $\sg$ by $i$. For
example, if $\sg = 2~7~5~4$, then $red(\sg) = 1~4~3~2$. Given a
permutation $\tau$ in the symmetric group $S_j$, we define a
permutation $\sg = \sg_1 \cdots \sg_n \in S_n$ to have a
$\tau$-match at place $i$ provided $red(\sg_i \cdots \sg_{i+j-1}) =
\tau$. (In the literature, $\tau$-match is also called an occurrence
of the consecutive pattern $\tau$.) Let $\tau\mbox{-}mch(\sigma)$ be
the number of $\tau$-matches in the permutation $\sg$. To prevent
confusion, we note that a permutation not having a $\tau$-match is
different than a permutation being $\tau$-avoiding.  A permutation
is called $\tau$-avoiding if there are no indices $i_1 < \cdots <
i_j$ such that $ red[\sg_{i_1} \cdots \sg_{i_j}] = \tau$. For
example, if $\tau = 2~1~4~3$, then the permutation $3~2~1~4~6~5$
does not have a $\tau$-match but it does not avoid $\tau$ since
$red[2~1~6~5] = \tau$.

In the case  where $|\tau|=2$, then $\tau\mbox{-}mch(\sigma)$
reduces to familiar permutation statistics. That is, if $\sg = \sg_1
\cdots \sg_n \in S_n$, let $Des(\sg) = \{i:\sg_i > \sg_{i+1}\}$ and
$Rise(\sg) = \{i:\sg_i < \sg_{i+1}\}$. Then it is easy to
see that $(2~1)\mbox{-}match(\sg) = des(\sg) = |Des(\sg)|$ and
$(1~2)\mbox{-}match(\sg) = rise(\sg) = |Rise(\sg)|$.

We can consider a more refined pattern matching condition where we
take into account conditions involving equivalence mod $k$ for some
integer $k \geq 2$. That is, suppose we fix $k \geq 2$ and we are
given some sequence of distinct integers $\tau = \tau_1 \cdots
\tau_j$. Then we say that a permutation $\sg = \sg_1 \cdots \sg_n
\in S_n$ has a $\tau$-$k$-equivalence match at place $i$ provided
$red(\sg_i \cdots \sg_{i+j-1}) = red(\tau)$ and for all $s \in \{0,
\ldots, j-1\}$, $\sg_{i+s} = \tau_{1+s} \mod k$. For example, if
$\tau = 1~2$ and $\sg = 5~1~7~4~3~6~8~2$, then $\sg$ has
$\tau$-matches starting at positions 2, 5, and 6. However, if $k=2$,
then only the $\tau$-match starting at position 5 is a
$\tau$-$2$-equivalence match. Let
$\tau\mbox{-}k\mbox{-}emch(\sigma)$ be the number of
$\tau$-$k$-equivalence matches in the permutation $\sg$.

More generally, if $\Upsilon$ is a set of sequences of distinct
integers of length $j$, then we say that a permutation $\sg = \sg_1
\cdots \sg_n \in S_n$ has a $\Upsilon$-$k$-equivalence match at
place $i$ provided there is a $\tau \in \Upsilon$ such that
$red(\sg_i \cdots \sg_{i+j-1}) = red(\tau)$ and for all $s \in \{0,
\ldots, j-1\}$, $\sg_{i+s} = \tau_{1+s} \mod k$. Let
$\Upsilon\mbox{-}k\mbox{-}emch(\sigma)$ be the number of
$\Upsilon$-$k$-equivalence matches in the permutation $\sg$ and
$\Upsilon\mbox{-}k\mbox{-}enlap(\sigma)$ be the maximum number of
non-overlapping $\Upsilon$-$k$-equivalence matches in $\sg$.

One can then study the following polynomials
\begin{eqnarray*}
T_{\tau,k,n}(x) &=& \sum_{\sg \in S_n}
x^{\tau\mbox{-}k\mbox{-}emch(\sigma)} =
\sum_{s=0}^n T_{\tau,k,n}^s x^s \ \mbox{and} \label{eq:k-tau} \\
U_{\Upsilon,k,n}(x) &=& \sum_{\sg \in S_n}
x^{\Upsilon\mbox{-}k\mbox{-}emch(\sigma)} =
\sum_{s=0}^n U_{\Upsilon,k,n}^s x^s.\label{eq:k-Up}
\end{eqnarray*}

Clearly our polynomials $A^{(k)}_n(x)$ and $B^{(k)}_n(x)$ are special cases of $U_{\Upsilon,k,n}(x)$ where all the patterns have length 2.
Papers \cite{liese} and \cite{liese-rem} have started the study of these more general types of polynomials in the case where all the patterns have length 2.


\begin{thebibliography}{4}

\bibitem{comtet} L. Comtet: {\em Advanced Combinatorics}, D. Reidel Publishing Co., Dordrecht, 1974.


\bibitem{Gasper} G. Gasper, Summation formulas for basic hypergeometric series,
{\em Siam. J. Math. Anal.}, {\bf 12} (1981), 196-200.

\bibitem{haglund2} James Haglund, Rook theory and hypergeometric Series,
{\em Advances in Applied Math.}, {\bf 17} (1996), 408-459.

\bibitem{haglund} James Haglund, private communication.



\bibitem{kitrem} S. Kitaev and J. Remmel: Classifying descents
according to parity, preprint.

\bibitem{liese} Jeffrey Liese, Classifying ascents and descents with
specified equivalences mod k,
Proceedings of 18-th International Conference on Formal Power Series and Algebra Combinatorics, San Diego, CA (2006).

\bibitem{liese-rem} Jeffrey Liese and Jeffrey B. Remmel, Pattern matching in permutaions relative to equivalence mod $k$, in preparation.

\bibitem{macmah} P. A. MacMahon: {\em Combinatory Analysis}, Vol. 1
and 2, Cambridge Univ. Press, Cambridge, 1915 (reprinted by
Chelsea, New York, 1955).

\bibitem{A=B} M. Petkovsek, H. S. Wilf and D. Zeilberger: $A=B$.
Wellesley, MA: A. K. Peters.

\end{thebibliography}
\end{document}